\documentclass[11pt,twoside,a4paper]{article}
\usepackage{amsmath}
\usepackage{amssymb}
\usepackage[latin1]{inputenc}
\usepackage[T1]{fontenc}   
\usepackage[english]{babel}
\newcommand{\tun}{\begin{picture}(5,0)(-2,-1)
\put(0,0){\circle*{2}}
\end{picture}}

\newcommand{\tdeux}{\begin{picture}(7,7)(0,-1)
\put(3,0){\circle*{2}}
\put(3,0){\line(0,1){5}}
\put(3,5){\circle*{2}}
\end{picture}}

\newcommand{\ttroisun}{\begin{picture}(15,8)(-5,-1)
\put(3,0){\circle*{2}}
\put(-0.65,0){$\vee$}
\put(6,7){\circle*{2}}
\put(0,7){\circle*{2}}
\end{picture}}
\newcommand{\ttroisdeux}{\begin{picture}(5,12)(-2,-1)
\put(0,0){\circle*{2}}
\put(0,0){\line(0,1){5}}
\put(0,5){\circle*{2}}
\put(0,5){\line(0,1){5}}
\put(0,10){\circle*{2}}
\end{picture}}

\newcommand{\tquatreun}{\begin{picture}(15,12)(-5,-1)
\put(3,0){\circle*{2}}
\put(-0.65,0){$\vee$}
\put(6,7){\circle*{2}}
\put(0,7){\circle*{2}}
\put(3,7){\circle*{2}}
\put(3,0){\line(0,1){7}}
\end{picture}}
\newcommand{\tquatredeux}{\begin{picture}(15,18)(-5,-1)
\put(3,0){\circle*{2}}
\put(-0.65,0){$\vee$}
\put(6,7){\circle*{2}}
\put(0,7){\circle*{2}}
\put(0,14){\circle*{2}}
\put(0,7){\line(0,1){7}}
\end{picture}}
\newcommand{\tquatretrois}{\begin{picture}(15,18)(-5,-1)
\put(3,0){\circle*{2}}
\put(-0.65,0){$\vee$}
\put(6,7){\circle*{2}}
\put(0,7){\circle*{2}}
\put(6,14){\circle*{2}}
\put(6,7){\line(0,1){7}}
\end{picture}}
\newcommand{\tquatrequatre}{\begin{picture}(15,18)(-5,-1)
\put(3,5){\circle*{2}}
\put(-0.65,5){$\vee$}
\put(6,12){\circle*{2}}
\put(0,12){\circle*{2}}
\put(3,0){\circle*{2}}
\put(3,0){\line(0,1){5}}
\end{picture}}
\newcommand{\tquatrecinq}{\begin{picture}(9,19)(-2,-1)
\put(0,0){\circle*{2}}
\put(0,0){\line(0,1){5}}
\put(0,5){\circle*{2}}
\put(0,5){\line(0,1){5}}
\put(0,10){\circle*{2}}
\put(0,10){\line(0,1){5}}
\put(0,15){\circle*{2}}
\end{picture}}

\newcommand{\tcinqun}{\begin{picture}(20,8)(-5,-1)
\put(3,0){\circle*{2}}
\put(-0.5,0){$\vee$}
\put(6,7){\circle*{2}}
\put(0,7){\circle*{2}}
\put(3,0){\line(2,1){10}}
\put(3,0){\line(-2,1){10}}
\put(-7,5){\circle*{2}}
\put(13,5){\circle*{2}}
\end{picture}}
\newcommand{\tcinqdeux}{\begin{picture}(15,14)(-5,-1)
\put(3,0){\circle*{2}}
\put(-0.65,0){$\vee$}
\put(6,7){\circle*{2}}
\put(0,7){\circle*{2}}
\put(3,7){\circle*{2}}
\put(3,0){\line(0,1){7}}
\put(0,7){\line(0,1){7}}
\put(0,14){\circle*{2}}
\end{picture}}
\newcommand{\tcinqtrois}{\begin{picture}(15,15)(-5,-1)
\put(3,0){\circle*{2}}
\put(-0.65,0){$\vee$}
\put(6,7){\circle*{2}}
\put(0,7){\circle*{2}}
\put(3,7){\circle*{2}}
\put(3,0){\line(0,1){7}}
\put(3,7){\line(0,1){7}}
\put(3,14){\circle*{2}}
\end{picture}}
\newcommand{\tcinqquatre}{\begin{picture}(15,14)(-5,-1)
\put(3,0){\circle*{2}}
\put(-0.65,0){$\vee$}
\put(6,7){\circle*{2}}
\put(0,7){\circle*{2}}
\put(3,7){\circle*{2}}
\put(3,0){\line(0,1){7}}
\put(6,7){\line(0,1){7}}
\put(6,14){\circle*{2}}
\end{picture}}
\newcommand{\tcinqcinq}{\begin{picture}(15,19)(-5,-1)
\put(3,0){\circle*{2}}
\put(-0.65,0){$\vee$}
\put(6,7){\circle*{2}}
\put(0,7){\circle*{2}}
\put(6,14){\circle*{2}}
\put(6,7){\line(0,1){7}}
\put(0,14){\circle*{2}}
\put(0,7){\line(0,1){7}}
\end{picture}}
\newcommand{\tcinqsix}{\begin{picture}(15,20)(-7,-1)
\put(3,0){\circle*{2}}
\put(-0.65,0){$\vee$}
\put(6,7){\circle*{2}}
\put(0,7){\circle*{2}}
\put(-3.65,7){$\vee$}
\put(3,14){\circle*{2}}
\put(-3,14){\circle*{2}}
\end{picture}}
\newcommand{\tcinqsept}{\begin{picture}(15,20)(-5,-1)
\put(3,0){\circle*{2}}
\put(-0.65,0){$\vee$}
\put(6,7){\circle*{2}}
\put(0,7){\circle*{2}}
\put(2.35,7){$\vee$}
\put(3,14){\circle*{2}}
\put(9,14){\circle*{2}}
\end{picture}}
\newcommand{\tcinqhuit}{\begin{picture}(15,26)(-5,-1)
\put(3,0){\circle*{2}}
\put(-0.65,0){$\vee$}
\put(6,7){\circle*{2}}
\put(0,7){\circle*{2}}
\put(0,14){\circle*{2}}
\put(0,7){\line(0,1){7}}
\put(0,21){\circle*{2}}
\put(0,14){\line(0,1){7}}
\end{picture}}
\newcommand{\tcinqneuf}{\begin{picture}(15,26)(-5,-1)
\put(3,0){\circle*{2}}
\put(-0.65,0){$\vee$}
\put(6,7){\circle*{2}}
\put(0,7){\circle*{2}}
\put(6,14){\circle*{2}}
\put(6,7){\line(0,1){7}}
\put(6,21){\circle*{2}}
\put(6,14){\line(0,1){7}}
\end{picture}}
\newcommand{\tcinqdix}{\begin{picture}(15,19)(-5,-1)
\put(3,5){\circle*{2}}
\put(-0.5,5){$\vee$}
\put(6,12){\circle*{2}}
\put(0,12){\circle*{2}}
\put(3,0){\circle*{2}}
\put(3,0){\line(0,1){12}}
\put(3,12){\circle*{2}}
\end{picture}}
\newcommand{\tcinqonze}{\begin{picture}(15,26)(-5,-1)
\put(3,5){\circle*{2}}
\put(-0.65,5){$\vee$}
\put(6,12){\circle*{2}}
\put(0,12){\circle*{2}}
\put(3,0){\circle*{2}}
\put(3,0){\line(0,1){5}}
\put(0,12){\line(0,1){7}}
\put(0,19){\circle*{2}}
\end{picture}}
\newcommand{\tcinqdouze}{\begin{picture}(15,26)(-5,-1)
\put(3,5){\circle*{2}}
\put(-0.65,5){$\vee$}
\put(6,12){\circle*{2}}
\put(0,12){\circle*{2}}
\put(3,0){\circle*{2}}
\put(3,0){\line(0,1){5}}
\put(6,12){\line(0,1){7}}
\put(6,19){\circle*{2}}
\end{picture}}
\newcommand{\tcinqtreize}{\begin{picture}(5,26)(-2,-1)
\put(0,0){\circle*{2}}
\put(0,0){\line(0,1){7}}
\put(0,7){\circle*{2}}
\put(0,7){\line(0,1){7}}
\put(0,14){\circle*{2}}
\put(-3.65,14){$\vee$}
\put(-3,21){\circle*{2}}
\put(3,21){\circle*{2}}
\end{picture}}
\newcommand{\tcinqquatorze}{\begin{picture}(9,26)(-5,-1)
\put(0,0){\circle*{2}}
\put(0,0){\line(0,1){5}}
\put(0,5){\circle*{2}}
\put(0,5){\line(0,1){5}}
\put(0,10){\circle*{2}}
\put(0,10){\line(0,1){5}}
\put(0,15){\circle*{2}}
\put(0,15){\line(0,1){5}}
\put(0,20){\circle*{2}}
\end{picture}}

\newcommand{\tdun}[1]{\begin{picture}(10,5)(-2,-1)
\put(0,0){\circle*{2}}
\put(3,-2){\tiny #1}
\end{picture}}

\newcommand{\bun}{
\begin{picture}(5,5)(-2,0)
\put(0,0){\line(0,0){5}}
\end{picture}}

\newcommand{\bdeux}{
\begin{picture}(10,15)(-5,0)
\put(-3.5,4.1){$\vee$}
\put(0,0){\line(0,1){5}}
\end{picture}}

\newcommand{\btroisun}{
\begin{picture}(15,20)(-7,0)
\put(-3.5,4.1){$\vee$}
\put(0,0){\line(0,1){5}}
\put(-6.3,10.7){$\vee$}
\end{picture}}
\newcommand{\btroisdeux}{
\begin{picture}(15,20)(-5,0)
\put(-3.5,4.1){$\vee$}
\put(0,0){\line(0,1){5}}
\put(-0.7,10.7){$\vee$}
\end{picture}}

\newcommand{\bquatreun}{
\begin{picture}(15,25)(-7,0)
\put(-3.5,4.1){$\vee$}
\put(0,0){\line(0,1){5}}
\put(-6.3,10.7){$\vee$}
\put(-9.1,17.3){$\vee$}
\end{picture}}
\newcommand{\bquatredeux}{
\begin{picture}(15,25)(-7,0)
\put(-3.5,4.1){$\vee$}
\put(0,0){\line(0,1){5}}
\put(-6.3,10.7){$\vee$}
\put(-3.5,17.3){$\vee$}
\end{picture}}
\newcommand{\bquatretrois}{
\begin{picture}(15,25)(-7,0)
\put(-3.5,4.1){$\vee$}
\put(0,0){\line(0,1){5}}
\put(-0.7,10.7){$\vee$}
\put(-3.5,17.3){$\vee$}
\end{picture}}
\newcommand{\bquatrequatre}{
\begin{picture}(15,25)(-7,0)
\put(-3.5,4.1){$\vee$}
\put(0,0){\line(0,1){5}}
\put(-0.7,10.7){$\vee$}
\put(2.1,17.3){$\vee$}
\end{picture}}
\newcommand{\bquatrecinq}{
\begin{picture}(15,20)(-7,0)
\put(0,0){\line(0,1){5}}
\put(0,5){\line(2,1){10}}
\put(0,5){\line(-2,1){10}}
\put(-3,6.5){\line(0,1){5}}
\put(3,6.5){\line(0,1){5}}
\end{picture}}

\newcommand{\bddeux}[1]{
\begin{picture}(25,22)(-12,0)
\put(0,0){\line(0,0){10}}
\put(0,10){\line(1,1){10}}
\put(0,10){\line(-1,1){10}}
\put(-9,7){\tiny #1}
\end{picture}}

\newcommand{\bdtroisun}[2]{
\begin{picture}(30,40)(-20,0)
\put(0,0){\line(0,0){10}}
\put(0,10){\line(1,1){10}}
\put(0,10){\line(-1,1){20}}
\put(-10,20){\line(1,1){10}}
\put(-9,7){\tiny #1}
\put(-19,17){\tiny #2}
\end{picture}}
\newcommand{\bdtroisdeux}[2]{
\begin{picture}(30,40)(-10,0)
\put(0,0){\line(0,0){10}}
\put(0,10){\line(1,1){20}}
\put(0,10){\line(-1,1){10}}
\put(10,20){\line(-1,1){10}}
\put(3,7){\tiny #1}
\put(13,17){\tiny #2}
\end{picture}}

\input{xy}
\xyoption{all}

\newcommand{\tdelta}{\tilde{\Delta}}
\newcommand{\h}{{\cal H}}
\newcommand{\m}{{\cal M}}
\newcommand{\F}{\mathbf{F}}
\newcommand{\T}{\mathbf{T}}
\newcommand{\TT}{\mathbb{T}}
\newcommand{\D}{{\cal D}}
\renewcommand{\P}{\mathbb{P}}
\newcommand{\Q}{\mathbb{PRIM}}
\newcommand{\Vect}{\mathbb{VECT}}
\newcommand{\haut}{\nearrow}
\newcommand{\bas}{\searrow}
\newcommand{\fleche}{{\stackrel{?}{\rightarrow}}}
\newcommand{\circhaut}{{\nearrow \hspace{-3.9mm} \circ \hspace{2mm}}}
\newcommand{\circbas}{{\searrow \hspace{-3.9mm} \circ \hspace{2mm}}}
\newcommand{\circfleche}{{\fleche \hspace{-3.0mm} \circ \hspace{2mm}}}
\newcommand{\prodhaut}{{\nearrow \hspace{-3.9mm} \bullet \hspace{2mm}}}
\newcommand{\prodbas}{{\searrow \hspace{-3.9mm} \bullet \hspace{2mm}}}
\newcommand{\prodfleche}{{\fleche \hspace{-3.0mm} \bullet \hspace{2mm}}}

\title{The infinitesimal Hopf algebra and the operads of planar forests}
\date{}
\author{Loïc Foissy \\ \\
{\small{\it Laboratoire de Mathématiques, FRE3111, Université de Reims}}\\
\small{{\it Moulin de la Housse - BP 1039 - 51687 REIMS Cedex 2, France}}\\
\small{e-mail : loic.foissy@univ-reims.fr}}

\newtheorem{defi}{\indent Definition}
\newtheorem{lemma}[defi]{\indent Lemma}
\newtheorem{cor}[defi]{\indent Corollary}
\newtheorem{theo}[defi]{\indent Theorem}
\newtheorem{prop}[defi]{\indent Proposition}

\newenvironment{proof}{{\bf Proof.}}{\hfill $\Box$}

\begin{document}

\maketitle

ABSTRACT. We introduce two operads which own the set of planar forests as a basis. With its usual product and two other products 
defined by different types of graftings, the algebra of planar rooted trees $\h$ becomes an algebra over these operads.
The compatibility with the infinitesimal coproduct of $\h$ and these structures is studied. As an application, an inductive way of computing 
the dual basis of $\h$ for its infinitesimal pairing is given.  Moreover, three Cartier-Quillen-Milnor-Moore theorems are given for the operads 
of planar forests and a rigidity theorem for one of them.\\


AMS CLASSIFICATION. 16W30, 05C05, 18D50.

\tableofcontents

\section*{Introduction}

The Connes-Kreimer Hopf algebra of rooted trees, introduced in \cite{Connes,Kreimer1,Kreimer2,Kreimer3},
is a commutative, non cocommutative Hopf algebra, its coproduct being given by admissible cuts of trees.
A non commutative version, the Hopf algebra of planar rooted trees, is introduced in \cite{Foissy2,Holtkamp}.
We furthemore introduce in \cite{FoissyInf} an infinitesimal version of this object, replacing admissible cuts by left admissible cuts:
this last object is here denoted by $\h$. Similarly with the Hopf case, $\h$ is a self-dual object and it owns a non-degenerate, symmetric
Hopf pairing, denoted by $\langle-,-\rangle$. This pairing is related to a partial order on the set of planar rooted forests,
making it isomorphic to the Tamari poset. As a consequence, $\h$ is given a dual basis denoted by $(f_F)_{F\in \F}$,
indexed by the set $\F$ of planar forest. In particular, the sub-family $(f_t)_{t\in \T}$ indexed by the set of planar rooted trees
$\T$ is a basis of the space of primitive elements of $\h$.\\

The aim of this text is to introduce two structures of operad on the space of planar forests.
We introduce two (non-symmetric) operads $\P_\bas$ and $\P_\haut$ defined in the following way:
\begin{enumerate}
\item $\P_\bas$ is generated by $m$ and $\bas \in \P_\bas(2)$, with relations:
$$ \left\{ \begin{array}{rcl}
m\circ (\bas,I)&=&\bas \circ (I,m),\\
m\circ (m,I)&=&m\circ (I,m),\\
\bas \circ (m,I)&=&\bas \circ (I,\bas).
\end{array}\right.$$
\item $\P_\haut$ is generated by $m$ and $\haut \in \P_\haut(2)$, with relations:
$$ \left\{ \begin{array}{rcl}
m\circ (\haut,I)&=&\haut \circ (I,m),\\
m\circ (m,I)&=&m\circ (I,m),\\
\haut \circ (\haut,I)&=&\haut \circ (I,\haut).
\end{array}\right.$$
\end{enumerate}
We then introduce two products on $\h$ or on its augmentation ideal $\m$, denoted by $\haut$ and $\bas$. The product $F \haut G$
consists of grafting $F$ on the left leave of $G$ and the product $F \bas G$ consists of grafting $F$ on the left root of $G$.
Together with its usual product $m$, $\m$ becomes both a $\P_\bas$- and a $\P_\haut$-algebra. More precisely,  $\m$ is the free $\P_\bas$-  
and $\P_\haut$-algebra generated by a single element $\tun$. As a consequence, $\P_\bas$ and $\P_\haut$ inherits a combinatorial
representation using planar forests, with composition iteratively described using the products $\bas$ and $\haut$.

We then give several applications of these operadic structures. For example, the antipode of $\h$ is described in term of the operad $\P_\bas$. 
We show how to compute elements $f_t$'s, with $t\in \T$, using the action of $\P_\bas$, and the elements $f_F$'s, $F\in \F$ from the preceding ones
using  the action of $\P_\haut$. Combining all these results, it is possible to compute by induction the basis $(f_F)_{F\in \F}$.

We finally study the compatibilities of products $m$, $\haut$, $\bas$, the coproduct $\tdelta$ and the coproduct $\tdelta_\haut$ dual of $\haut$.
This leads to the definition of two types of $\P_\haut$-bialgebras, and one type of $\P_\bas$-bialgebras.
Each type then define a suboperad of $\P_\haut$ or $\P_\bas$ corresponding to primitive elements of $\m$, which are explicitively described:
\begin{enumerate}
\item The first one is a free operad, generated by the element $\tdeux-\tun\tun\in \P_\haut(2)$. As a consequence, the
space of primitive elements of $\h$ admits a basis $(p_t)_{t\in \T_b}$ indexed by the set of planar binary trees.
The link with the basis $(f_t)_{t\in \T}$ is given with the help of the Tamari order.
\item The second one admits a combinatorial representation in terms of planar rooted trees. 
It is generated by the corollas $c_n\in \P_\haut(n)$, $n\geq 2$, with the following relations: for all $k,l \geq 2$,
$$c_k \circ (c_l,\underbrace{I,\ldots,I}_{\mbox{\scriptsize $k-1$ times}})=c_l \circ (\underbrace{I,\ldots,I}_{\mbox{\scriptsize $l-1$ times}},c_k).$$
\item The third one admits a combinatorial representation in terms of planar rooted trees, and is freely generated by $\tdeux\in \P_\bas(2)$.
\end{enumerate}

We also give the definition of a double $\P_\haut$-bialgebra, combining the two types of $\P_\haut$-bialgebras already introduced.
We then prove a rigidity theorem: any double $\P_\haut$-bialgebra connected as a coalgebra is isomorphic to a decorated version of $\m$.\\

This text is organised as follows: the first section gives several recalls on the infinitesimal Hopf algebra of planar rooted trees and its pairing. 
The two products $\bas$ and $\haut$ are introduced in section 2, as well as the combinatorial representation of the two associated operads. 
The applications to the computation of $(f_F)_{F\in \F}$ is given in section 3. Section 4 is devoted to the study of the suboperads 
of primitive elements and the last section deals with the rigidity theorem for double $\P_\haut$-bialgebras. \\

{\bf Notations.} \begin{enumerate}
\item We shall denote by $K$ a commutative field, of any characteristic. Every vector space, algebra, coalgebra, etc, will be taken over $K$.
\item Let $(A,\Delta,\varepsilon)$ be a counitary coalgebra. Let $1 \in A$, non zero, such that $\Delta(1)=1\otimes 1$.
We then define the non counitary coproduct:
$$\tdelta : \left\{ \begin{array}{rcl}
Ker(\varepsilon)&\longrightarrow & Ker(\varepsilon) \otimes Ker(\varepsilon)\\
a&\longrightarrow & \tdelta(a)=\Delta(a)-a\otimes 1-1\otimes a.
\end{array} \right.$$ 
We shall use the Sweedler notations $\Delta(a)=a^{(1)}\otimes a^{(2)}$ and $\tdelta(a)=a'\otimes a''$.
\end{enumerate}

\section{Planar rooted forests and their infinitesimal Hopf algebra}

We here recall some results and notations of \cite{FoissyInf}.

\subsection{Planar trees and forests}

\begin{enumerate}
\item The set of planar trees is denoted by $\T$, and the set of planar forests is denoted by $\F$.
The weight of a planar forest is the number of its vertices. For all $n \in \mathbb{N}$, we denote by $\F(n)$ the set of planar forests of weight $n$.

{\bf Examples.}
Planar rooted trees of weight $\leq 5$:
$$\tun,\tdeux,\ttroisun,\ttroisdeux,\tquatreun, \tquatredeux,\tquatretrois,\tquatrequatre,\tquatrecinq,\tcinqun,\tcinqdeux,\tcinqtrois,\tcinqquatre,\tcinqcinq,
\tcinqsix,\tcinqsept,\tcinqhuit,\tcinqneuf,\tcinqdix,\tcinqonze,\tcinqdouze,\tcinqtreize,\tcinqquatorze.$$
Planar rooted forests of weight $\leq 4$:
$$1,\tun,\tun\tun,\tdeux,\tun\tun\tun,\tdeux\tun,\tun \tdeux,\ttroisun,\ttroisdeux,\tun\tun\tun\tun,\tdeux\tun\tun,\tun \tdeux \tun, \tun \tun \tdeux,
\ttroisun\tun,\tun \ttroisun,\ttroisdeux\tun,\tun \ttroisdeux,\tdeux\tdeux,\tquatreun,\tquatredeux,\tquatretrois,\tquatrequatre,\tquatrecinq.$$

\item The algebra $\h$ is the free associative, unitary algebra generated by $\T$. As a consequence, a linear basis of $\h$ is given by $\F$,
and its product is given by the concatenation of planar forests. 
\item We shall also need two partial orders and a total order on the set $Vert(F)$ of vertices of $F\in \F$, defined in \cite{Foissy2,FoissyInf}.
We put $F=t_1\ldots t_n$, and let $s,s'$ be two vertices of $F$.
\begin{enumerate}
\item  We shall say that $s \geq_{high} s'$ if there exists a path from $s'$ to $s$ in $F$, the edges of $F$ being oriented from the roots to the leaves. 
Note that $\geq_{high}$ is a partial order, whose Hasse graph is the forest $F$.
\item If $s$ and $s'$ are not comparable for $\geq_{high}$, we shall say that $s \geq_{left} s'$ if one of these assertions is satisfied:
\begin{enumerate}
\item $s$ is a vertex of $t_i$ and $s'$ is a vertex of $t_j$, with $i<j$.
\item $s$ and $s'$ are vertices of the same $t_i$, and $s \geq_{left} s'$ in the forest obtained from $t_i$ by deleting its root.
\end{enumerate}
This defines the partial order $\geq_{left}$ for all forests $F$, by induction on the the weight.
\item We shall say that $s \geq_{h,l} s'$ if $s \geq_{high} s'$ or $s \geq_{left} s'$. This defines a total order on the vertices of $F$. 
\end{enumerate}
\end{enumerate}

\subsection{Infinitesimal Hopf algebra of planar forests}

\begin{enumerate}
\item Let $F \in \F$. An {\it admissible cut} is a non empty cut of certain edges and trees of $F$, such that each path in a non-cut tree of $F$ 
meets at most one cut edge. The set of admissible cuts of $F$ will be denoted by $Adm(F)$. If $c$ is an admissible cut of $F$, 
the forest of the vertices which are over the cuts of $c$ will be denoted by $P^c(t)$ (branch of the cut $c$), and the remaining forest 
will be denoted by $R^c(t)$ (trunk of the cut). An admissible cut of $F$ will be said to be {\it left-admissible} if, for all vertices $x$ and $y$ of $F$,
$x\in P^c(F)$ and $x\leq_{left} y$ imply that $y\in P^c(F)$. The set of left-admissible cuts of $F$ will be denoted by $Adm^l(F)$.
\item $\h$ is given a coproduct by the following formula: for all $F\in \F$,
$$\Delta(F)=\sum_{c \in {\cal A}dm^l(F)} P^c(F) \otimes R^c(F)+F \otimes 1+1 \otimes F.$$
Then $(\h,\Delta)$ is an infinitesimal bialgebra, that is to say: for all $x,y\in \h$,
$$\Delta(xy)=(x\otimes 1)\Delta(y)+\Delta(x)(1\otimes y)-x\otimes y.$$

{\bf Examples.} \begin{eqnarray*}
\Delta(\tun)&=&\tun\otimes 1+1\otimes \tun,\\
\Delta(\tun\tun)&=&\tun\tun\otimes 1+1\otimes \tun\tun+\tun\otimes \tun,\\
\Delta(\tdeux)&=& \tdeux\otimes 1+1\otimes \tdeux+\tun \otimes \tun,\\
\Delta(\tdeux\tun)&=& \tdeux\tun\otimes 1+1\otimes \tdeux\tun+\tun \otimes \tun\tun+\tdeux\otimes \tun,\\
\Delta(\ttroisun)&=&\ttroisun\otimes 1+1\otimes \ttroisun+\tun\tun \otimes \tun+ \tun\otimes \tdeux,\\
\Delta(\ttroisdeux)&=&\ttroisdeux\otimes 1+1\otimes \ttroisdeux+\tdeux \otimes \tun+ \tun\otimes \tdeux,\\
\Delta(\tun\tun\tun\tun)&=&\tun\tun\tun\tun\otimes 1+1\otimes \tun\otimes \tun\tun\tun+\tun\tun\otimes \tun\tun+\tun\tun\tun\otimes \tun,\\
\Delta(\tdeux\tun\tun)&=&\tdeux\tun\tun\otimes 1+1\otimes \tdeux\tun\tun+\tun\otimes \tun\tun\tun+\tdeux\otimes \tun\tun+\tdeux\tun\otimes\tun,\\
\Delta(\tun\tdeux\tun)&=&\tun\tdeux\tun\otimes 1+1\otimes \tun\tdeux\tun+\tun\otimes \tdeux\tun+\tun\tun\otimes \tun\tun+\tun\tdeux\otimes\tun,\\
\Delta(\tun\tun\tdeux)&=&\tun\tun\tdeux\otimes 1+1\otimes \tun\tun\tdeux+\tun\otimes\tun \tdeux+\tun\tun\otimes \tdeux+\tun\tun\tun\otimes\tun,\\
\Delta(\tun\ttroisun)&=&\tun\ttroisun\otimes 1+1\otimes \tun\ttroisun+\tun\otimes \ttroisun+\tun\tun \otimes \tdeux+\tun\tun\tun\otimes \tun,\\
\Delta(\tun\ttroisdeux)&=&\tun\ttroisdeux\otimes 1+1\otimes \tun\ttroisdeux+\tun\otimes \ttroisdeux+\tun\tun \otimes \tdeux+\tun\tdeux\otimes \tun,\\
\Delta(\ttroisun\tun)&=&\ttroisun\tun\otimes 1+1\otimes \ttroisun+\tun\otimes\tdeux\tun+\tun\tun\otimes \tun\tun+\ttroisun\otimes \tun,\\
\Delta(\ttroisdeux\tun)&=&\ttroisdeux\tun\otimes 1+1\otimes \ttroisdeux+\tun\otimes\tdeux\tun+\tdeux\otimes \tun\tun+\ttroisdeux\otimes \tun,\\
\Delta(\tdeux\tdeux)&=&\tdeux\tdeux\otimes 1+1\otimes \tdeux\tdeux+\tun \otimes \tun\tdeux+\tdeux\otimes\tdeux+\tdeux\tun\otimes \tun,
\end{eqnarray*}
\begin{eqnarray*}
\Delta(\tquatreun)&=&\tquatreun\otimes 1+1\otimes\tquatreun+\tun\otimes \ttroisun+\tun\tun\otimes \tdeux+\tun\tun\tun\otimes \tun,\\
\Delta(\tquatredeux)&=&\tquatredeux\otimes 1+1\otimes\tquatredeux+\tun\otimes \ttroisun+\tdeux\otimes \tdeux+\tdeux\tun\otimes \tun,\\
\Delta(\tquatretrois)&=&\tquatretrois\otimes 1+1\otimes\tquatretrois+\tun\otimes \ttroisdeux+\tun\tun\otimes \tdeux+\tun\tdeux\otimes \tun,\\
\Delta(\tquatrequatre)&=&\tquatrequatre\otimes 1+1\otimes \tquatrequatre\tun\otimes \ttroisdeux+\tun\tun\otimes \tdeux+\ttroisun\otimes \tun,\\
\Delta(\tquatrecinq)&=&\tquatrecinq\otimes1 +1\otimes \tquatrecinq+\tun\otimes \ttroisdeux+\tdeux\otimes\tdeux+\ttroisdeux\otimes \tun.
\end{eqnarray*} \end{enumerate}

We proved in \cite{FoissyInf} that $\h$ is an infinitesimal Hopf algebra, that is to say has an antipode $S$. This antipode satisfies
$S(1)=1$, $S(t)  \in Prim(\h)$ for all $t\in \T$, and $S(F)=0$ for all $F \in \F-(\T\cup \{1\})$.

\subsection{Pairing on $\h$}

\begin{enumerate}
\item We define the operator $B^+:\h\longrightarrow \h$, which associates, to a forest $F\in \F$, the tree obtained by grafting the roots 
of the trees of $F$ on a common root. For example, $B^+(\ttroisun\tun)=\tcinqsix$, and $B^+(\tun\ttroisun)=\tcinqsept$.
\item The application $\gamma$ is defined by:
$$\gamma: \left\{ \begin{array}{rcl}
\h&\longrightarrow & \h\\
t_1\ldots t_n \in \F&\longrightarrow & \delta_{t_1,\tun} t_2\ldots t_n.
\end{array}\right.$$
\item There exists a unique pairing $\langle-,-\rangle:\h \times \h \longrightarrow K$, satisfying:
\begin{description}
\item[\textnormal{i.}] $\langle1,x\rangle =\varepsilon(x)$ for all $x \in \h$.
\item[\textnormal{ii.}] $\langle xy,z\rangle =\langle y\otimes x, \Delta(z)\rangle $ for all $x,y,z \in \h$.
\item[\textnormal{iii.}] $\langle B^+(x),y\rangle =\langle x,\gamma(y)\rangle $ for all $x,y \in \h$.
\end{description}
Moreover:
\begin{description}
\item[\textnormal{iv.}] $\langle-,-\rangle $ is symmetric and non-degenerate.
\item[\textnormal{v.}] If $x$ and $y$ are homogeneous of different weights, $\langle x,y\rangle =0$.
\item[\textnormal{vi.}] $\langle S(x),y\rangle =\langle x,S(y)\rangle $ for all $x,y \in \h$.
\end{description}
This pairing admits a combinatorial interpretation using the partial orders $\geq_{left}$ and $\geq_{high}$ and is related to the Tamari order on
 planar binary trees, see \cite{FoissyInf}.
\item We denote by $(f_F)_{F\in \F}$ the dual basis of the basis of forests for the pairing $\langle-,-\rangle$. In other terms, for all $F \in \F$, 
$f_F$ is defined by $\langle f_F,G\rangle =\delta_{F,G}$, for all forest $G \in \F$. The family $(f_t)_{t\in \T}$ is a basis of the space $Prim(\h)$
of primitive elements of $\h$.
\end{enumerate}

\section{The operads of forests and graftings}

\subsection{A few recalls on non-$\Sigma$-operads}

\begin{enumerate}
\item We shall work here with non-$\Sigma$-operads \cite{Markl}. Recall that such an object is a family $\P=(\P(n))_{n\in \mathbb{N}}$ 
of vector spaces, together with a composition for all $n,k_1,\ldots,k_n \in \mathbb{N}$:
$$\left\{ \begin{array}{rcl}
\P(n)\otimes \P(k_1)\otimes \ldots \otimes \P(k_n)&\longrightarrow & \P(k_1+\ldots+k_n)\\
p\otimes p_1\otimes \ldots \otimes p_n&\longrightarrow &p\circ(p_1,\ldots,p_n).
\end{array} \right.$$
The following associativity condition is satisfied: for all $p\in \P(n)$, $p_1 \in \P(k_1)$, $\ldots$, $p_n \in \P(k_n)$,  $p_{1,1},\ldots,p_{n,k_n} \in \P$,
\begin{eqnarray*}
&&(p\circ(p_1,\ldots,p_n))\circ (p_{1,1},\ldots,p_{1,k_1},\ldots,p_{n,1},\ldots,p_{n,k_n})\\
&=&p\circ (p_1\circ(p_{1,1},\ldots,p_{1,k_1}),\ldots,p_n\circ(p_{n,1},\ldots,p_{n,k_n})).
\end{eqnarray*}
Moreover, there exists a unit element $I\in \P(1)$, satisfying: for all $p\in \P(n)$,
$$\left\{ \begin{array}{rcl}
p\circ(I,\ldots,I)&=&p,\\
I\circ p&=&p.
\end{array}\right.$$
An operad is a non-$\Sigma$-operad $\P$ with a right action of the symmetric group $S_n$ on $\P(n)$ for all $n$,
satisfying a certain compatibility with the composition.
\item Let $\P$ be a non-$\Sigma$-operad. A $\P$-algebra is a vector space $A$, together with an action of $\P$:
$$ \left\{ \begin{array}{rcl}
\P(n)\otimes A^{\otimes n}&\longrightarrow & A\\
p\otimes a_1\otimes \ldots \otimes a_n&\longrightarrow &p.(a_1,\ldots,a_n),
\end{array} \right.$$
satisfying the following compatibility: for all $p\in \P(n)$, $p_1 \in \P(k_1)$, $\ldots$, $p_n \in \P(k_n)$, for all $a_{1,1},\ldots,a_{n,k_n} \in A$,
\begin{eqnarray*}
&&(p\circ(p_1,\ldots,p_n)).(a_{1,1},\ldots,a_{1,k_1},\ldots,a_{n,1}\ldots,a_{n,k_n})\\
&=&p.(p_1.(a_{1,1},\ldots,a_{1,k_1}),\ldots,p_n.(a_{n,1},\ldots,a_{n,k_n})).
\end{eqnarray*}
Moreover, $I.a=a$ for all $a\in A$.

In particular, if $V$ is a vector space, the free $\P$-algebra generated by $V$ is:
$$F_\P(V)=\bigoplus_{n\in \mathbb{N}} \P(n) \otimes V^{\otimes n},$$
with the action of $\P$ given by:
\begin{eqnarray*}
&&p.\left((p_1\otimes a_{1,1}\otimes \ldots \otimes a_{1,k_1}),\ldots,(p_n\otimes a_{n,1}\otimes \ldots \otimes a_{n,k_n})\right)\\
&=&(p\circ(p_1,\ldots,p_n)) \otimes a_{1,1}\otimes \ldots \otimes a_{1,k_1}\otimes \ldots \otimes a_{n,1}\otimes \ldots \otimes a_{n,k_n}.
\end{eqnarray*}

\item Let $\T_b$ be the set of planar binary trees:
$$\T_b=\left\{\bun,\bdeux,\btroisun, \btroisdeux, \bquatreun,\bquatredeux,\bquatretrois,\bquatrequatre,\bquatrecinq\ldots\right\}.$$
For all $n \in \mathbb{N}$, $\TT_b(n)$ is the vector space generated by the elements of $\T_b$ with $n$ leaves:
\begin{eqnarray*}
\TT_b(0)&=&(0),\\
\TT_b(1)&=&Vect\left(\bun\right),\\
\TT_b(2)&=&Vect\left(\bdeux\right),\\
\TT_b(3)&=&Vect\left(\btroisun, \btroisdeux \right),\\
\TT_b(4)&=&Vect\left(\:\bquatreun,\bquatredeux,\bquatretrois,\bquatrequatre,\bquatrecinq \right).\\
\end{eqnarray*}

The family of vector spaces $\TT_b$ is given a structure of non-$\Sigma$-operad by graftings on the leaves. More precisely, 
if $t,t_1,\ldots,t_n \in \T_b$, $t$ with $n$ leaves, then $t\circ (t_1,\ldots,t_n)$ is the binary tree obtained by grafting $t_1$ on the first leave of $t$,
$t_2$ on the second leave of $t$, and so on (note that the leaves of $t$ are ordered from left to right). The unit is $\bun$.

It is known that $\TT_b$ is the free non-$\Sigma$-operad generated by $\bdeux \in \TT_b(2)$. Similarly, given elements 
$m_1,\ldots,m_k$ in $\P(2)$, it is possible to describe the free non-$\Sigma$-operad $\P$ generated by these elements in terms of 
planar binary trees whose internal vertices are decorated by $m_1,\ldots,m_k$.
\end{enumerate}

\subsection{Presentations of the operads of forests}

\begin{defi}\textnormal{ \begin{enumerate}
\item $\P_\bas$ is the non-$\Sigma$-operad generated by $m$ and $\bas \in \P_\bas(2)$, with relations:
$$ \left\{ \begin{array}{rcl}
m\circ (\bas,I)&=&\bas \circ (I,m),\\
m\circ (m,I)&=&m\circ (I,m),\\
\bas \circ (m,I)&=&\bas \circ (I,\bas).
\end{array}\right.$$
\item $\P_\haut$ is the non-$\Sigma$-operad generated by $m$ and $\haut \in \P_\haut(2)$, with relations:
$$ \left\{ \begin{array}{rcl}
m\circ (\haut,I)&=&\haut \circ (I,m),\\
m\circ (m,I)&=&m\circ (I,m),\\
\haut \circ (\haut,I)&=&\haut \circ (I,\haut).
\end{array}\right.$$
\end{enumerate}} \end{defi}

{\bf Remark.} We shall prove in \cite{FoissyOperades} that these quadratic operads are Koszul.

\subsection{Grafting on the root}

Let $F,G \in \F-\{1\}$. We put $G=t_1\ldots t_n$ and $t_1=B^+(G_1)$. We define:
$$F\bas G=B^+(FG_1)t_2\ldots t_n.$$
In other terms, $F$ is grafted on the root of the first tree of $G$, on the left. In particular, $F \bas \tun=B^+(F)$. \\

{\bf Examples.} 
$$\begin{array}{rclcl|rclcl|rclcl|rclcl}
\tun\tun\tun &\bas& \tdeux&=&\tcinqun&\tdeux &\bas& \tun \tun \tun&=&\ttroisdeux \tun \tun&
\tun\tun&\bas&\tun\tun\tun&=&\ttroisun\tun\tun&\tun\tun\tun&\bas&\tun\tun&=&\tquatreun\tun\\
\tun\tdeux &\bas& \tdeux&=&\tcinqtrois&\tdeux &\bas& \tun \tdeux&=&\ttroisdeux \tdeux&
\tun\tun&\bas&\tun\tdeux&=&\ttroisun\tdeux&\tun\tdeux&\bas&\tun\tun&=&\tquatretrois\tun\\
\tdeux\tun &\bas& \tdeux&=&\tcinqdeux&\tdeux &\bas& \tdeux\tun &=&\tquatredeux \tun&
\tun\tun&\bas&\tdeux\tun&=&\tquatreun\tun& \tdeux\tun&\bas&\tun\tun&=&\tquatredeux\tun\\
\ttroisun &\bas& \tdeux&=&\tcinqsix&\tdeux &\bas& \ttroisun&=&\tcinqdeux&
\tun\tun&\bas&\ttroisun&=&\tcinqun&\ttroisun&\bas&\tun\tun&=&\tquatrequatre\tun\\
\ttroisdeux &\bas& \tdeux&=&\tcinqhuit&\tdeux &\bas& \ttroisdeux&=&\tcinqcinq&
\tun\tun&\bas& \ttroisdeux&=&\tcinqquatre& \ttroisdeux&\bas&\tun\tun&=&\tquatrecinq\tun.
\end{array}$$

Obviously, $\bas$ can be inductively defined in the following way:
for $F,G,H \in \F-\{1\}$,
$$\left\{ \begin{array}{rcl}
F\bas \tun&=&B^+(F),\\
F\bas (GH)&=& (F\bas G)H\\
F\bas B^+(G)&=&B^+(FG).
\end{array}\right.$$\\

We denote by $\m$ the augmentation ideal of $\h$, that is to say the vector space generated by the elements of $\F-\{1\}$.
We extend $\bas:\m\otimes \m \longrightarrow \m$ by linearity.

\begin{prop} \label{2}
For all $x,y,z \in \m$:
\begin{eqnarray}
\label{E1} x \bas (yz)&=&(x\bas y)z,\\
\label{E2} x\bas (y\bas z)&=&(xy) \bas z.
\end{eqnarray} \end{prop}

\begin{proof} We can restrict ourselves to $x,y,z \in \F-\{1\}$. Then (\ref{E1}) is immediate. In order to prove (\ref{E2}), we put $z=B^+(z_1)z_2$, 
$z_1,z_2\in \F$. Then:
$$x\bas (y\bas z)=x \bas (B^+(yz_1)z_2)=B^+(xyz_1)z_2=(xy) \bas (B^+(z_1)z_2)=(xy)\bas z,$$
which proves (\ref{E2}). \end{proof}

\begin{cor}
$\m$ is given a graded $\P_\bas$-algebra structure by its products $m$ and by $\bas$.
\end{cor}

\begin{proof} Immediate, by proposition \ref{2}. \end{proof}

\subsection{Grafting on the left leave}

Let $F,G \in \F$. Suppose that $G \neq 1$. Then $F \haut G$ is the planar forest obtained by grafting $F$ on the leave of $G$ 
which is at most on the left. For $G=1$, we put $F \haut 1=F$. In particular, $F \haut \tun=B^+(F)$.\\

{\bf Examples.}
$$\begin{array}{rclcl|rclcl|rclcl|rclcl}
\tun\tun\tun &\haut& \tdeux&=&\tcinqdix&\tdeux &\haut& \tun \tun \tun&=&\ttroisdeux \tun \tun&
\tun\tun&\haut&\tun\tun\tun&=&\ttroisun\tun\tun&\tun\tun\tun&\haut&\tun\tun&=&\tquatreun\tun\\
\tun\tdeux &\haut& \tdeux&=&\tcinqdouze&\tdeux &\haut& \tun \tdeux&=&\ttroisdeux \tdeux&
\tun\tun&\haut&\tun\tdeux&=&\ttroisun\tdeux&\tun\tdeux&\haut&\tun\tun&=&\tquatretrois\tun\\
\tdeux\tun &\haut& \tdeux&=&\tcinqonze&\tdeux &\haut& \tdeux\tun &=&\tquatrecinq \tun&
\tun\tun&\haut&\tdeux\tun&=&\tquatrequatre\tun& \tdeux\tun&\haut&\tun\tun&=&\tquatredeux\tun\\
\ttroisun &\haut& \tdeux&=&\tcinqtreize&\tdeux &\haut& \ttroisun&=&\tcinqhuit&
\tun\tun&\haut&\ttroisun&=&\tcinqsix&\ttroisun&\haut&\tun\tun&=&\tquatrequatre\tun\\
\ttroisdeux &\haut& \tdeux&=&\tcinqquatorze&\tdeux &\haut& \ttroisdeux&=&\tcinqquatorze&
\tun\tun&\haut& \ttroisdeux&=&\tcinqtreize& \ttroisdeux&\haut&\tun\tun&=&\tquatrecinq\tun.
\end{array}$$

In an obvious way, $\haut$ can be inductively defined in the following way: for $F,G,H \in \F$,
$$\left\{\begin{array}{rcl}
F\haut 1&=&F,\\
F\haut (GH)&=& (F\haut G)H\mbox{ if } G \neq 1,\\
F\haut B^+(G)&=&B^+(F\haut G).
\end{array} \right.$$
We extend $\haut:\h\otimes \h \longrightarrow \h$ by linearity.

\begin{prop} \label{4} \begin{enumerate}
\item For all $x,z \in \h$, $y \in \m$:
\begin{equation}
\label{E3} x \haut (yz)=(x\haut y)z.
\end{equation}
\item For all $x,y,z \in \h$:
$$x\haut (y\haut z)=(x\haut y) \haut z.$$
So $(\h, \haut)$ is an associative algebra, with unitary element $1$.
\end{enumerate} \end{prop}

\begin{proof} Note that (\ref{E3}) is immediate for $x,y,z \in \F$, with $y\neq 1$. This implies the first point.
In order to prove the second point, we consider:
$$Z=\{z\in \h\:/\: \forall x,y \in \h, \: x\haut (y\haut z)=(x\haut y) \haut z \}.$$
Let us first prove that $1 \in Z$: for all $x,y \in \h$,
$$x \haut (y \haut 1)=x \haut y=(x\haut y) \haut 1.$$
Let $z_1,z_2 \in Z$. Let us show that $z_1 z_2 \in Z$. By linearity, we can separate the proof into two cases:
\begin{enumerate}
\item $z_1=1$. Then it is obvious.
\item $\varepsilon(z_1)=0$. Let $x,y \in \h$. By the first point: 
\begin{eqnarray*}
x \haut( y\haut (z_1z_2))&=&x \haut ((y\haut z_1)z_2))\\
&=&(x \haut (y\haut z_1))z_2\\
&=&((x \haut y) \haut z_1)z_2\\
&=&(x\haut y) \haut (z_1z_2). 
\end{eqnarray*} \end{enumerate}
So $Z$ is a subalgebra of $\h$. Let us show that it is stable by $B^+$. Let $z\in Z$, $x,y \in \h$. Then:
\begin{eqnarray*}
x \haut (y\haut B^+(z))&=&x \haut B^+(y\haut z)\\
&=&B^+ (x \haut (y\haut z))\\
&=&B^+ ((x \haut y)\haut z)\\
&=&(x \haut y) \haut B^+ (z).
\end{eqnarray*}
So $Z$ is a subalgebra of $\h$, stable by $B^+$. Hence, $Z=\h$. \end{proof}\\

{\bf Remarks.} \begin{enumerate}
\item (\ref{E3}) is equivalent to: for any $x,y,z\in \h$, 
$$x\haut(yz)-\varepsilon(y) x \haut z=(x\haut y)z-\varepsilon(y)xz.$$
\item Let $F\in \F-\{1\}$. There exists a unique family $(\tun F_1,\ldots,\tun F_n)$ of elements of $\F$ such that: 
$$F=(\tun F_1)\haut\ldots\haut(\tun F_n).$$
For example, $\tcinqsix\tdeux \tun=(\tun \tun)\haut (\tun\tun)\haut (\tun \tdeux \tun)$.
As a consequence, $(\h,\haut)$ is freely generated by $\tun \F$ as an associative algebra.
\end{enumerate}

\begin{cor}
$\m$ is given a graded $\P_\haut$-algebra structure by its product $m$ and by $\haut$.
\end{cor}

\begin{proof} Immediate, by proposition \ref{4}. \end{proof}

\subsection{Dimensions of $\P_\bas$ and $\P_\haut$}

We now compute the dimensions of $\P_\bas(n)$ and $\P_\haut(n)$ for all $n$ and deduce that $\m$ is the free $\P_\bas$- and $\P_\haut$-algebra
generated by $\tun$.\\

{\bf Notation.} We denote by $r_n$ the number of planar rooted forests and we put $\displaystyle R(X)=\sum_{n=1}^{+\infty} r_n X^n$.
It is well-known (see \cite{Foissy2,Stanley2}) that $R(X)=\displaystyle \frac{1-2X-\sqrt{1-4X}}{2X}$.

\begin{prop}
For $\fleche \in \{\bas,\haut\}$ and all $n\in \mathbb{N}^*$, in the $\P_\fleche$-algebra $\m$:
$$\P_\fleche(n).(\tun,\ldots,\tun)=Vect(\mbox{planar forests of weight $n$}).$$
As a consequence, $\m$ is generated as a $\P_\fleche$-algebra by $\tun$.
\end{prop}

\begin{proof} $\subseteq$. Immediate, as $\m$ is a graded $\P_\fleche$-algebra. \\

$\supseteq$. Induction on $n$. For $n=1$, $I.(\tun)=\tun$. For $n\geq 2$, two cases are possible.
\begin{enumerate}
\item $F=F_1F_2$, $weight(F_i)=n_i<n$. By the induction hypothesis, there exists $p_1,p_2 \in \P_\fleche$, such that $F_1=p_1.(\tun, \ldots,\tun)$ 
and $F_2=p_2.(\tun,\ldots,\tun)$. Then $(m\circ(p_1,p_2)).(\tun,\ldots,\tun)=m.(F_1,F_2)=F_1F_2$.
\item $F\in \T$. Let us put $F=B^+(G)$. Then there exists $p\in \P_\fleche$, such that $p.(\tun,\ldots,\tun)=G$. Then:
$$\left\{\begin{array}{rcccl}
(\bas \circ (p,I)).(\tun,\ldots,\tun)&=&G\bas \tun&=&F,\\
(\haut \circ (p,I)).(\tun,\ldots,\tun)&=&G\haut \tun&=&F.
\end{array}\right.$$
\end{enumerate}
Hence, in both cases, $F \in \P_\fleche(n).(\tun,\ldots,\tun)$. \end{proof}

\begin{cor}
For all $\fleche \in \{\bas,\haut\}$, $n \in \mathbb{N}^*$, $dim(\P_\fleche(n))\geq r_n$.
\end{cor}

\begin{proof} Because we proved the surjectivity of the following application:
$$ ev_\fleche: \left\{\begin{array}{rcl}
\P_\fleche (n)& \longrightarrow &Vect(\mbox{planar forests of weight $n$})\\
p & \longrightarrow &p.(\tun, \ldots,\tun).
\end{array}\right. $$ \end{proof}

\begin{lemma}
For all $\fleche \in \{\bas,\haut\}$, $n \in \mathbb{N}^*$, $dim(\P_\fleche(n))\leq r_n$.
\end{lemma}

\begin{proof} We prove it for $\fleche=\haut$. Let us fix $n \in \mathbb{N}^*$. Then $\P_\haut(n)$ is linearly generated by planar binary trees 
whose internal vertices are decorated by $m$ and $\haut$. The following relations hold:
$$\begin{array}{rclcrclcrcl}
\bdtroisun{$\haut $}{$\haut $}&=&\bdtroisdeux{$\haut $}{$\haut $}&,&
\bdtroisun{$m $}{$\haut $}&=&\bdtroisdeux{$\haut $}{$m $}&,&
\bdtroisun{$m $}{$m $}&=&\bdtroisdeux{$m $}{$m $}.
\end{array}$$
In the sequel of the proof, we shall say that such a tree is {\it admissible} if it satisfies the following conditions:
\begin{enumerate}
\item For each internal vertex $s$ decorated by $m$, the left child of $s$ is a leave.
\item For each internal vertex $s$ decorated by $\haut$, the left child of $s$ is a leave or is decorated by $m$.
\end{enumerate}
For example, here are the admissible trees with $1$, $2$ or $3$ leaves:
$$\bun, \bddeux{$m$}, \bddeux{$\haut $}, \bdtroisun{$\haut $}{$m$}, \bdtroisdeux{$m$}{$m$}, \bdtroisdeux{$\haut $}{$m$}, 
\bdtroisdeux{$m$}{$\haut $}, \bdtroisdeux{$\haut $}{$\haut $}.$$

The preceding relations imply that $\P_\haut(n)$ is linearly generated by admissible trees with $n$ leaves. 
So $dim(\P_\haut(n))$ is smaller than $a_n$, the number of admissible trees with $n$ leaves. 
For $n \geq 2$, we denote by $b_n$ the number of admissible trees with $n$ leaves whose root is decorated by $m$,
 and by $c_n$ the number of admissible trees with $n$ leaves whose root is decorated by $\haut$.
We also put $b_1=1$ and $c_1=0$. Finally, we define:
$$ A(X)=\sum_{n\geq 1} a_n X^n,\hspace{1cm}  B(X)=\sum_{n\geq 1} b_n X^n,\hspace{1cm} C(X)=\sum_{n\geq 1} c_n X^n.$$ 

Immediately, $A(X)=B(X)+C(X)$. Every admissible tree with $n\geq 2$ leaves whose root is decorated by $m$ is of the form $m \circ (I,t)$, 
where $t$ is an admissible tree with $n-1$ leaves. Hence, $B(X)=XA(X)+X$. Moreover, every admissible tree with $n\geq 2$ leaves 
whose root is decorated by $\haut$ is of the form $\haut \circ (t_1,t_2)$, where $t_1$ is an admissible tree with $k$ leaves whose
eventual root is decorated by $m$ and $t_2$ an admissible tree with $n-k$ leaves ($1\leq k\leq n-1$). 
Hence, for all $n \geq 2$, $\displaystyle c_n=\sum_{k=1}^{n-1} b_k a_{n-k}$, so $C(X)=B(X)A(X)$. As a conclusion:
$$\left\{\begin{array}{rcl}
A(X)&=& B(X)+C(X),\\
B(X)&=&XA(X)+X,\\
C(X)&=&B(X)A(X).
\end{array}\right. $$
So $A(X)=XA(X)+X+B(X)A(X)=XA(X)+X+XA(X)^2+XA(X)$, and:
$$XA(X)^2+(2X-1)A(X)+X=0.$$ 
As $a_1=1$: 
$$ A(X)=\frac{1-2X-\sqrt{1-4X}}{2X}=R(X).$$
So, for all $n \geq 1$, $dim(\P_\haut(n))\leq a_n=r_n$. The proof is similar for $\P_\bas$. \end{proof}\\

As immediate consequences:

\begin{theo} \label{9}
For $\fleche \in \{\bas,\haut\}$, $n \in \mathbb{N}^*$, $dim(\P_\fleche(n))=r_n$. Moreover, the following application is bijective:
$$ ev_\fleche: \left\{\begin{array}{rcl}
\P_\fleche (n)& \longrightarrow &Vect(\mbox{planar forests of weight $n$})\subseteq \m\\
p & \longrightarrow &p.(\tun, \ldots,\tun).
\end{array}\right. $$
\end{theo}

\begin{cor} \begin{enumerate}
\item $(\m,m,\bas)$ is the free $\P_\bas$-algebra generated by $\tun$.
\item $(\m,m,\haut)$ is the free $\P_\haut$-algebra generated by $\tun$.
\end{enumerate} \end{cor}

\subsection{A combinatorial description of the composition}

Let $\fleche \in \{\bas,\haut\}$. We identify $\P_\fleche$ and the vector space of non-empty planar forests via theorem \ref{9}. 
In other terms, we identify $F \in \F(n)$ and $ev_\fleche^{-1}(F) \in \P_\fleche(n)$.

{\bf Notations.} \begin{enumerate}
\item In order to distinguish the compositions in $\P_\bas$ and $\P_\haut$, we now denote:
\begin{enumerate}
\item $F \circbas(F_1,\ldots,F_n)$ the composition of $\P_\bas$,
\item $F \circhaut(F_1,\ldots,F_n)$ the composition of $\P_\haut$.
\end{enumerate}
\item In order to distinguish the action of the operads  $\P_\bas$ and $\P_\haut$ on $\m$, we now denote:
\begin{enumerate}
\item $F\prodbas(x_1,\ldots,x_n)$ the action of $\P_\bas$ on $\m$,
\item $F\prodhaut(x_1,\ldots,x_n)$ the action of $\P_\haut$ on $\m$.
\end{enumerate} \end{enumerate}

Our aim in this paragraph is to describe the compositions of $\P_\bas$ and $\P_\haut$ in term of forests. We shall prove the following result:

\begin{theo} \label{11} \begin{enumerate}
\item The composition of $\P_\bas$ in the basis of planar forests can be inductively defined in this way:
$$  \left\{ \begin{array}{rcl}
\tun\circbas (H)&=&H,\\
B^+(F)\circbas (H_1,\ldots, H_{n+1})
&=&(F\circbas (H_1,\ldots,H_n))\bas H_{n+1},\\
FG\circbas (H_1,\ldots, H_{n_1+n_2})&=&F\circbas(H_1,\ldots,H_{n_1})G\circbas (H_{n_1+1},\ldots, H_{n_1+n_2}).
\end{array} \right.$$
\item The composition of $\P_\haut$ in the basis of planar forests can be inductively defined in this way:
$$  \left\{ \begin{array}{rcl}
\tun\circhaut (H)&=&H,\\
B^+(F)\circhaut (H_1,\ldots, H_{n+1})
&=&(F\circhaut (H_1,\ldots,H_n))\haut H_{n+1},\\
FG\circhaut (H_1,\ldots, H_{n_1+n_2})&=&F\circhaut(H_1,\ldots,H_{n_1})G\circhaut (H_{n_1+1},\ldots, H_{n_1+n_2}).
\end{array} \right.$$
\end{enumerate} \end{theo}

{\bf Examples.} Let $F_1,F_2,F_3 \in \F-\{1\}$.
$$\begin{array}{rclc|crcl}
\tun\tun\circhaut (F_1,F_2)&=&F_1F_2,&&&\tun\tun\circbas (F_1,F_2)&=&F_1F_2,\\
\tdeux\circhaut (F_1,F_2)&=&F_1\haut F_2,&&&\tdeux\circbas (F_1,F_2)&=&F_1\bas F_2,\\
\tun\tun\tun\circhaut (F_1,F_2,F_3)&=&F_1F_2F_3,&&&\tun\tun\tun\circbas (F_1,F_2,F_3)&=&F_1F_2F_3,\\
\tun\tdeux\circhaut (F_1,F_2,F_3)&=&F_1(F_2\haut F_3),&&&\tun\tdeux\circbas (F_1,F_2,F_3)&=&F_1(F_2\bas F_3),\\
\tdeux\tun\circhaut (F_1,F_2,F_3)&=&(F_1\haut F_2)F_3,&&&\tdeux\tun\circbas (F_1,F_2,F_3)&=&(F_1\bas F_2)F_3,\\
\ttroisun\circhaut (F_1,F_2,F_3)&=&(F_1F_2)\haut F_3,&&&\ttroisun\circbas (F_1,F_2,F_3)&=&(F_1F_2)\bas F_3,\\
\ttroisdeux\circhaut (F_1,F_2,F_3)&=&(F_1\haut F_2)\haut F_3,&&&\ttroisdeux\circbas (F_1,F_2,F_3)&=&(F_1\bas F_2)\bas F_3.
\end{array}$$

\begin{prop} \label{12}
Let $\fleche \in \{\bas,\haut\}$.
\begin{enumerate}
\item  $\tun$ is the unit element of $\P_\fleche$.
\item $\tun \tun=m$ in $\P_\fleche(2)$. Consequently, in $\P_\fleche$, $\tun\tun \circ (F,G)=FG$ for all $F,G \in \F-\{1\}$.
\item Let $F,G \in \F$. In $\P_\fleche$, $\tdeux=\fleche$. Consequently, $\tdeux \circfleche (F,G)=F\fleche G$ for all $F,G \in \F-\{1\}$.
\end{enumerate} \end{prop}

\begin{proof} \begin{enumerate}
\item Indeed, $ev_\fleche(\tun)=\tun=ev_\fleche(I)$. Hence, $\tun=I$. 
\item By definition, $ev_\fleche(\tun \tun)=\tun \tun=ev_\fleche(m)$. So $\tun \tun=m$ in $\P_\fleche(2)$. Moreover, for all $F,G \in \F-\{1\}$:
\begin{eqnarray*}
ev_\fleche(FG)&=&FG\\
&=&m\prodfleche (F,G)\\
&=&m\prodfleche(F\prodfleche(\tun,\ldots,\tun),G\prodfleche(\tun,\ldots,\tun))\\
&=&\left(m\circfleche(F,G)\right)\prodfleche(\tun,\ldots,\tun)\\
&=&ev_\fleche(m\circfleche(F,G)).
\end{eqnarray*}
So $FG=m\circfleche(F,G)=\tun\tun\circfleche(F,G)$.
\item Indeed, $ev_\fleche(\tdeux)=\tun \fleche \tun=ev_\fleche(\fleche)$. So $\tdeux=\fleche$ in $\P_\fleche(2)$. Moreover:
\begin{eqnarray*}
ev_\fleche(F\fleche G)&=&F\fleche G\\
&=&\fleche\prodfleche(F,G)\\
&=&\fleche\prodfleche(F\prodfleche(\tun,\ldots,\tun),G\prodfleche(\tun,\ldots,\tun))\\
&=&(\fleche \circfleche(F,G)).(\tun,\ldots,\tun)\\
&=&ev_\fleche(\fleche \circfleche(F,G)).
\end{eqnarray*}
So, $F\fleche G=\fleche\circfleche(F,G)=\tdeux\circfleche(F,G)$. 
\end{enumerate} \end{proof}

\begin{prop} \label{13} \begin{enumerate}
\item Let $F,G\in \F$, different from $1$, of respective weights $n_1$ and $n_2$. Let $H_{1,1},\ldots,H_{1,n_1}$ 
and $H_{2,1},\ldots,H_{2,n_2}\in \F-\{1\}$. Let $\fleche\in \{\bas,\haut\}$. Then, in $\P_\fleche$:
$$(FG)\circfleche (H_{1,1},\ldots, H_{1,n_1},H_{2,1},\ldots,H_{2,n_2})=F\circfleche (H_{1,1},\ldots,H_{1,n_1})
G\circfleche (H_{2,1},\ldots,H_{2,n_2}).$$
\item Let $F\in \F$, of weight $n\geq 1$. Let $H_1,\ldots,H_{n+1}\in \F$. In $\P_\fleche$:
$$B^+(F)\circfleche(H_1,\ldots, H_{n+1})=(F\circfleche(H_1,\ldots,H_n))\fleche H_{n+1}.$$
\end{enumerate} \end{prop}

\begin{proof} \begin{enumerate}
\item Indeed, in $\P_\fleche$:
\begin{eqnarray*}
&&(FG) \circfleche (H_{1,1},\ldots, H_{1,n_1},H_{2,1},\ldots,H_{2,n_2})\\
&=&(m\circfleche (F,G))\circfleche (H_{1,1},\ldots, H_{1,n_1},H_{2,1},\ldots,H_{2,n_2})\\
&=&m\circfleche (F\circfleche (H_{1,1},\ldots, H_{1,n_1}), G\circfleche (H_{2,1},\ldots,H_{2,n_2}))\\
&=&F\circfleche (H_{1,1},\ldots, H_{1,n_1}) G\circfleche (H_{2,1},\ldots,H_{2,n_2})).
\end{eqnarray*}
\item In $\P_\fleche$:
\begin{eqnarray*}
B^+(F)\circfleche(H_1,\ldots,H_{n+1})&=&(F \fleche \tun)\circfleche(H_1,\ldots,H_{n+1})\\
&=&(\tdeux\circfleche(F,\tun))\circfleche(H_1,\ldots,H_{n+1})\\
&=&\tdeux\circfleche(F \circfleche (H_1,\ldots,H_n),\tun\circfleche(H_{n+1}))\\
&=&\tdeux\circfleche(F \circfleche (H_1,\ldots,H_n),H_{n+1})\\
&=&(F\circfleche(H_1,\ldots,H_n))\fleche H_{n+1}.
\end{eqnarray*}
 \end{enumerate} \end{proof}

Combining propositions \ref{12} and \ref{13}, we obtain theorem \ref{11}.

\section{Applications to the infinitesimal Hopf algebra $\h$}

\subsection{Antipode of $\h$}

We here give a description of the antipode of $\h$ in terms of the action $\prodbas$ of the operad $\P_\bas$.\\

{\bf Notations.} For all $n \in \mathbb{N}^*$, we denote $l_n=(B^+)^n(1)\in \F(n)$. For example:
$$l_1=\tun\:,\: l_2=\tdeux \:,\: l_3=\ttroisdeux \:,\: l_4=\tquatrecinq \:,\: l_5=\tcinqquatorze\ldots$$ 

\begin{lemma}
Let $t\in \T$. There exists a unique $k\in \mathbb{N}^*$, and a unique family $(t_2\ldots,t_k)\in \T^{k-1}$ such that:
$$t=l_k\prodbas(\tun,t_2,\ldots,t_k).$$
\end{lemma}

\begin{proof} Induction on the weight $n$ of $t$. If $n=1$, then $t=\tun$, so $k=1$ and the family is empty. We suppose the result at all rank $<n$. 
We put $t=B^+(s_1\ldots s_m)$. Necessarily, $t_k=B^+(s_2\ldots s_m)$ and $l_{n-1}\prodbas(\tun,t_2,\ldots,t_{k-1})=s_1$.
We conclude with the induction hypothesis on $s_1$. \end{proof}\\

{\bf Example.}
$$ \begin{picture}(10,26)(-2,-1)
\put(3,0){\circle*{2}}
\put(-.6,0){$\vee $}
\put(6,7){\circle*{2}}
\put(0,7){\circle*{2}}
\put(3,7){\circle*{2}}
\put(3,0){\line(0,1){7}}
\put(-3.65,7){$\vee $}
\put(-3.4,14){\circle*{2}}
\put(2.5,14){\circle*{2}}
\put(-7,14){$\vee $}
\put(-6.5,21){\circle*{2}}
\put(-0.5,21){\circle*{2}}
\put(-0.5,21){\line(0,1){7}}
\put(-0.5,28){\circle*{2}}
\end{picture}=l_4\prodbas(\tun,\ttroisdeux,\tdeux,\ttroisun).$$

\begin{defi}
\textnormal{For all $n \in \mathbb{N}^*$, we put 
$\displaystyle p_n = \sum_{k=1}^{n}\: \sum_{\substack{a_1+\ldots+a_k=n \\ \forall i, \:a_i> 0}}  (-1)^k l_{a_1}\ldots l_{a_k}$.}
 \end{defi}

{\bf Examples.} \begin{eqnarray*}
p_1&=&\tun,\\
p_2&=&-\tdeux+\tun\tun,\\
p_3&=&-\ttroisdeux+\tdeux\tun+\tun\tdeux-\tun\tun\tun,\\
p_4&=&-\tquatrecinq+\ttroisdeux\tun+\tdeux\tdeux+\tun \ttroisdeux-\tdeux\tun\tun-\tun\tdeux\tun-\tun\tun\tdeux+\tun\tun\tun\tun.
\end{eqnarray*}
Remark that $p_n$ is in fact the antipode of $l_n$ in $\h$. It is also the antipode of $l_n$ 
in the non commutative Connes-Kreimer Hopf algebra of planar trees \cite{Foissy2}.

\begin{cor}
Let $t \in \T$, written under the form $t=l_k\prodbas(t_1,\ldots,t_k)$, with $t_1=\tun$. Then:
$$S(t)= p_k\prodbas(t_1,\ldots,t_k).$$
\end{cor}

\begin{proof} Corollary of proposition 15 of \cite{FoissyInf}, observing that left cuts are cuts on edges from the root of $t_i$ to the root of $t_{i+1}$ in $t$,
for $i=1,\ldots, n-1$. \end{proof}

\subsection{Inverse of the application $\gamma$}

\label{s3-2}

\begin{prop}
The restriction $\gamma:Prim(\h)\longrightarrow \h$ is bijective.
\end{prop}

\begin{proof} By proposition 21 of \cite{FoissyInf}:
$$\gamma_{\mid Prim(\h)} : \left\{ \begin{array}{rcl}
Prim(\h)&\longrightarrow & \h\\
f_{B^+(F)}\:(F\in \F)&\longrightarrow &f_F.
\end{array}\right.$$
So this restriction is clearly bijective. \end{proof}\\

We shall denote $\gamma^{-1}_{\mid Prim(\h)}:\h \longrightarrow Prim(\h)$ the inverse of this restriction. Then, for all $F \in \F$, 
$\gamma^{-1}_{\mid Prim(\h)}(f_F)=f_{B^+(F)}$. Our aim is to express $\gamma^{-1}_{\mid Prim(\h)}$ in the basis of forests.\\

We define inductively a sequence $(q_n)_{n\in \mathbb{N}^*}$ of elements of $\P_\bas$:
$$  \left\{\begin{array}{rcl}
q_1&=&\tun \in \P_\bas(1),\\
q_2&=&\tun\tun-\tdeux \in \P_\bas(2),\\
q_{n+1}&=&(\tun\tun-\tdeux)\circbas(q_n,\tun) \in \P_\bas(n+1)\mbox{ for }n\geq 1.
\end{array}\right. $$
For all  $F\in \F$, $\tun\tun \circbas (F,\tun)=F\tun$ and $\tdeux\circbas(F,\tun)=B^+(F)$. So, $q_n$ can also be defined in the following way:
$$  \left\{\begin{array}{rcl}
q_1&=&\tun \in \P_\bas(1),\\
q_{n+1}&=&q_n \tun-B^+(q_n)\in \P_\bas(n+1)\mbox{ for } n\geq 1.
\end{array}\right. $$\\

{\bf Examples.}
\begin{eqnarray*}
q_3&=&\tun \tun \tun -\tdeux\tun -\ttroisun+\ttroisdeux,\\
q_4&=&\tun\tun \tun \tun -\tdeux\tun\tun -\ttroisun\tun+\ttroisdeux\tun-\tquatreun+\tquatredeux+\tquatrequatre-\tquatrecinq,\\
q_5&=&\tun\tun \tun \tun \tun-\tdeux\tun\tun\tun-\ttroisun\tun\tun+\ttroisdeux\tun\tun-\tquatreun\tun+\tquatredeux\tun+\tquatrequatre\tun-\tquatrecinq\tun \\
&&-\tcinqun+\tcinqdeux+\tcinqsix-\tcinqhuit+\tcinqdix-\tcinqonze-\tcinqtreize+\tcinqquatorze.
\end{eqnarray*}

\begin{lemma} \label{18}
Let $F\in \F-\{1\}$, and $t\in \T$. Then, in $\h$:
$$\Delta(F\bas t)=(F\bas t)\otimes 1+1\otimes (F\bas t)+F'\otimes F''\bas t+ Ft'\otimes t''+ F\otimes t.$$
\end{lemma}

\begin{proof}  The non-empty and non-total left-admissible cuts of the tree $F\bas t$ are:
\begin{description}
\item[\textnormal{-}] The cut on the edges relating $F$ to $t$. For this cut $c$, $P^c(F \bas t)=F$ and $R^c(F \bas t)=t$.
\item[\textnormal{-}] Cuts acting only on edges of $F$ or on edges relating $F$ to $t$, at the exception of the preceding case. For such a cut, 
there exists a unique non-empty, non-total left-admissible cut $c'$ of $F$,  such that $P^c(F \bas t)=P^{c'}(F)$ and $R^c(F \bas t)=R^{c'}(F) \bas t$.
\item[\textnormal{-}] Cuts acting on edges of $t$. Then necessarily $F \subseteq P^c(F \bas t)$. For such a cut, there exists a unique non-empty,
non-total left-admissible cut $c'$ of $t$, such that $P^c(F \bas t)=FP^{c'}(t)$ and $R^c(F \bas t)=R^{c'}(t)$.
\end{description}
Summing these cuts, we obtain the announced compatibility. \end{proof}

\begin{prop} \label{19} 
Let $F=t_1\ldots t_n \in \F$. Then:
$$\gamma^{-1}_{\mid Prim(\h)}(F)=q_{n+1}\prodbas(\tun,t_1,\ldots,t_n).$$
\end{prop}

\begin{proof} {\it First step.} Let us show the following property: for all $x \in Prim(\h)$, $t \in \T$, $q_2\prodbas(x,t)$ is primitive.
By lemma \ref{18}, using the linearity in $F$:
\begin{eqnarray*}
\Delta(x\bas t)&=&(x\bas t)\otimes 1+1\otimes (x\bas t)+x \otimes t+xt'\otimes t'',\\
\Delta(x t) &=&xt \otimes 1 +1 \otimes xt+x \otimes t+xt'\otimes t'',\\
\Delta(q_2\prodbas(x,t))&=&\Delta(xt-x\bas t)\\
&=&(xt-x\bas t) \otimes 1+1\otimes (xt-x\bas t).
\end{eqnarray*}

{\it Second step.} Let us show that for all $x \in Prim(\h)$, $t_1,\ldots,t_n\in \T$, $q_{n+1}\prodbas(x,t_1,\ldots,t_n)\in Prim(\h)$ by induction on $n$.
This is obvious for $n=0$, as $q_1\prodbas (x)=x$. Suppose the result at rank $n-1$. Then:
\begin{eqnarray*}
q_{n+1}\prodbas(x,t_1,\ldots,t_n) &=& (q_2\circbas (q_n,I))\prodbas (x,t_1,\ldots,t_n)\\
&=&q_2\prodbas(\underbrace{q_n\prodbas(x,t_1,\ldots,t_{n-1})}_{\in Prim(\h)},t_n)\:\in \: Prim(\h),
\end{eqnarray*}
by the first step. As the tree $\tun$ is primitive, we deduce that, for all forest $F=t_1\ldots t_n \in \F$, $q_{n+1}\prodbas(\tun,t_1,\ldots,t_n)\in Prim(\h)$.\\

{\it Third step.} Let us show that for all $x,y\in \m$, $\gamma(q_2\prodbas(x,y))=\gamma(x) y$. We can limit ourselves to $x,y\in \F-\{1\}$.
Then $q_2\prodbas(x,y)=xy-x\bas y$. Moreover, by definition of $\bas$, $x \bas y$ is a forest whose first tree is not equal to $\tun$.
Hence, $\gamma(q_2\prodbas(x,y))=\gamma(xy)-0=\gamma(x)y.$\\

{\it Last step.} Let us show by induction on $n$ that $\gamma(q_{n+1}\prodbas(\tun,t_1,\ldots,t_n)) =t_1\ldots t_n$. 
As $q_1\prodbas(\tun)=\tun$, this is obvious if $n=0$. Let us suppose the result at rank $n-1$. By the third step:
\begin{eqnarray*}
\gamma(q_{n+1}\prodbas(\tun,t_1,\ldots,t_n))
&=&\gamma(q_2\prodbas(q_n\prodbas (\tun,t_1,\ldots,t_{n-1}),t_n))\\
&=&\gamma(q_n\prodbas(\tun,t_1,\ldots,t_{n-1}))t_n\\
&=&t_1\ldots t_n.
\end{eqnarray*}
Consequently, $x=q_{n+1}\prodbas(\tun,t_1,\ldots,t_n)\in Prim(\h)$, and satisfies $\gamma(x)=t_1\ldots t_n$, which proves proposition \ref{19}. \end{proof}\\

{\bf Examples.} Let $t_1,t_2,t_3 \in \T$.
\begin{eqnarray*}
\gamma^{-1}_{\mid Prim(\h)}(t_1)&=&\tun t_1-\tun\bas t_1,\\
\gamma^{-1}_{\mid Prim(\h)}(t_1t_2)&=&\tun t_1t_2-(\tun\bas t_1)t_2-(\tun t_1)\bas t_2+(\tun\bas t_1)\bas t_2,\\
\gamma^{-1}_{\mid Prim(\h)}(t_1t_2t_3)&=&\tun t_1t_2t_3-(\tun\bas t_1)t_2t_3-(\tun t_1)\bas t_2t_3+(\tun\bas t_1)\bas t_2t_3-(\tun t_1t_2)\bas t_3\\
&&+(\tun\bas t_1t_2)\bas t_3+((\tun t_1)\bas t_2)\bas t_3-((\tun\bas t_1)\bas t_2)\bas t_3.
\end{eqnarray*}

\subsection{Elements of the dual basis}

\begin{lemma} \label{20}
For all $x,y\in \h$, $\Delta(x\haut y)=x\haut y^{(1)} \otimes y^{(2)} + x^{(1)} \otimes x^{(2)}\haut y-x\otimes y$. 
In other terms, $(\h,\haut,\Delta)$ is an infinitesimal Hopf algebra.
\end{lemma}

\begin{proof}  We restrict to $x=F\in \F-\{1\}$, $y=G \in \F-\{1\}$. The non-empty and non-total left-admissible cuts of the tree $F\haut G$ are:
\begin{description}
\item[\textnormal{-}] The cut on the edges relating $F$ to $G$. For this cut $c$, $P^c(F \haut G)=F$ and $R^c(F \haut G)=G$.
\item[\textnormal{-}] Cuts acting only on edges of $F$ or on edges relating $F$ to $G$, at the exception of the preceding case. For such a cut, 
there exists a unique non-empty, non-total left-admissible cut $c'$ of $F$,  such that $P^c(F \haut G)=P^{c'}(F)$  and $R^c(F \haut G)=R^{c'}(F) \haut G$.
\item[\textnormal{-}] Cuts acting on edges of $G$. Then necessarily $F \subseteq P^c(F \haut G)$. For such a cut, there exists a unique non-empty,
non-total left-admissible cut $c'$ of $t$, such that $P^c(F \haut G)=F\haut P^{c'}(G)$ and $R^c(F \haut G)=R^{c'}(G)$.
\end{description}
Summing these cuts, we obtain, denoting $\Delta(F)=F\otimes 1+1\otimes F+F'\otimes F''$ and $\Delta(G)=G\otimes 1+1\otimes G+G'\otimes G''$:
\begin{eqnarray*}
\tdelta(F \haut G)&=&(F \haut G) \otimes 1+1\otimes (F \haut G) +F\otimes G+F'\otimes F''\haut G+F \haut G'\otimes G''\\
&=&(F \otimes 1) \haut \Delta(G)+\Delta(F) \haut (1 \otimes G)-F \otimes G.
\end{eqnarray*}
So $(\h,\haut,\Delta)$ is an infinitesimal bialgebra. As it is graded and connected, it has an antipode. \end{proof}

\begin{prop} \label{21} 
Let $F=t_1\ldots t_n \in \F$. Then $f_F=f_{t_n}\haut \ldots \haut f_{t_1}$.
\end{prop}

\begin{proof} {\it First step.} We show the following result: for all $F \in \F$, $t\in \T$, $f_F \haut f_t=f_{tF}$. We proceed by induction 
on the weight $n$ of $F$.  If $n=0$, then $F=1$ and the result is obvious. We now suppose that the result is true at all rank $<n$. 
Let be $G \in \F$, and let us prove that $\langle f_F \haut f_t,G\rangle =\delta_{tF,G}$. Three cases are possible.
\begin{enumerate}
\item $G=1$. Then $\langle f_F\haut f_t,G\rangle =\langle f_F\haut f_t,1\rangle =\varepsilon(f_F\haut f_t)=0=\delta_{tF,G}$.
\item $G=G_1G_2$, $G_i\neq 1$. Then, by lemma \ref{20}:
\begin{eqnarray*}
\langle f_F\haut f_t,G\rangle&=&\langle \Delta(f_F\haut f_t),G_2\otimes G_1\rangle \\ 
&=&\sum_{F_1F_2=F} \langle f_{F_2} \otimes f_{F_1}\haut f_t,G_2\otimes G_1\rangle \\ 
&&+\langle f_F \haut f_t \otimes 1+f_F \haut 1 \otimes f_t ,G_2\otimes G_1\rangle -\langle f_F \otimes f_t ,G_2\otimes G_1\rangle \\
&=& \sum_{\substack{F_1F_2=F,\\ weight(F_1)<n}} \langle f_{F_2} \otimes f_{F_1}\haut f_t,G_2\otimes G_1\rangle 
+\langle 1 \otimes f_F\haut f_t,G_2\otimes G_1\rangle \\
&&+\langle f_F \haut f_t \otimes 1, G_2\otimes G_1 \rangle+\langle f_F \otimes f_t ,G_2\otimes G_1\rangle
-\langle f_F \otimes f_t ,G_2\otimes G_1\rangle \\
&=&  \sum_{\substack{F_1F_2=F,\\ weight(F_1)<n}} \langle f_{F_2} \otimes f_{tF_1},G_2\otimes G_1\rangle\\
&=& \sum_{\substack{F_1F_2=F,\\ weight(F_1)<n}} \delta_{F_2,G_2}\delta_{tF_1,G_1}\\
&=&\delta_{tF,G}.
\end{eqnarray*}
\item $G=B^+(G_1)$. Note that  $f_F \haut f_t$ is a linear span of forests $H_1\haut H_2$, with $H_1$, $H_2\neq 1$.
By definition of $\haut$, the first tree of such a forest is not $\tun$. Hence, $\gamma(f_F\haut f_t)=0$ and:
$$\langle f_F \otimes f_t,G\rangle =\langle \gamma(f_F \otimes f_t), G_1\rangle =0=\delta_{tF,G},$$
as $tF \notin \T$ because $F \neq 1$.\\
\end{enumerate}

{\it Second step.} We now prove proposition \ref{21} by induction on $n$. It is obvious for $n=1$. Suppose the result at rank $n-1$. By the first step:
$$f_{t_1\ldots t_n}=f_{t_2\ldots t_n}\haut f_{t_1}=(f_{t_n} \haut \ldots \haut f_{t_2})\haut f_{t_1}= f_{t_n} \haut \ldots \haut f_{t_2}\haut f_{t_1},$$
using the induction hypothesis for the second equality. \end{proof}\\

{\bf Remarks.}
\begin{enumerate}
\item As an immediate corollary, because $\haut$ is associative, for all forests $F_1,\ldots,F_k\in \F$, 
$f_{F_1\ldots F_k}=f_{F_k} \haut \ldots \haut f_{F_1}$.
\item In term of operads, proposition \ref{21} can be rewritten in the following way:
\end{enumerate}

\begin{cor}
Let $F_1,\ldots,F_n \in \F$. Then $f_{F_1\ldots F_n}= l_n\circhaut(f_{F_n}, \ldots, f_{F_1})$.
\end{cor}

{\bf Remark.} Hence, the dual basis $(f_F)_{F\in \F}$ can be inductively computed, using proposition 21 of \cite{FoissyInf}, 
together with propositions \ref{19} and \ref{21} of the present text:
$$  \left\{ \begin{array}{rcl}
f_1&=&1,\\
f_{t_1\ldots t_n}&=&f_{t_n}\haut \ldots \haut f_{t_1},\\
f_{B^+(t_1\ldots t_n)}&=&\gamma^{-1}_{\mid Prim(\h)}(f_{t_1\ldots t_n}).
\end{array}\right.$$ 
For example:
$$\begin{array}{|c|}
\hline \begin{array}{rclc|rcl}
f_1&=&1&&
f_{\tun}&=&\tun\\
f_{\tun\tun}&=&\tdeux&&
f_{\tdeux}&=&-\tdeux+\tun\tun\\
f_{\tun\tun\tun}&=&\ttroisdeux&&
f_{\tun\tdeux}&=&-\ttroisdeux+\ttroisun\\
f_{\tdeux\tun}&=&-\ttroisdeux+\tdeux\tun&&
f_{\ttroisun}&=&-\ttroisun+\tun\tdeux\\
f_{\ttroisdeux}&=&\ttroisdeux-\tdeux\tun-\tun\tdeux+\tun\tun\tun&&
f_{\tun\tun\tun\tun}&=&\tquatrecinq\\
f_{\tun\tun\tdeux}&=&-\tquatrecinq+\tquatrequatre&&
f_{\tun\tdeux\tun}&=&-\tquatrecinq+\tquatredeux\\
f_{\tun\ttroisun}&=&-\tquatrequatre+\tquatretrois&&
f_{\tun\ttroisdeux}&=&\tquatrecinq-\tquatredeux-\tquatretrois+\tquatreun\\
f_{\tdeux\tun\tun}&=&-\tquatrecinq+\ttroisdeux\tun&&
f_{\tdeux\tdeux}&=&\tquatrecinq-\tquatrequatre-\ttroisdeux\tun+\ttroisun\tun\\
f_{\ttroisun\tun}&=&-\tquatredeux+\tdeux\tdeux&&
f_{\ttroisdeux\tun}&=&\tquatrecinq-\ttroisdeux\tun-\tdeux\tdeux+\tdeux\tun\tun\\
f_{\tquatreun}&=&-\tquatretrois+\tun\ttroisdeux&&
f_{\tquatretrois}&=&\tquatretrois-\tquatreun-\tun\ttroisdeux+\tun\ttroisun\\
f_{\tquatredeux}&=&\tquatrequatre-\ttroisun\tun-\tun\ttroisdeux+\tun\tdeux\tun&&
f_{\tquatrequatre}&=&\tquatredeux-\tdeux\tdeux-\tun\ttroisun+\tun\tun\tdeux\\
\end{array}\\
\hline f_{\tquatrecinq}=-\tquatrecinq+\ttroisdeux\tun+\tdeux\tdeux-\tdeux\tun\tun+
\tun\ttroisdeux-\tun\tdeux\tun-\tun\tun\tdeux+\tun\tun\tun\tun.\\
\hline \end{array}$$

\section{Primitive suboperads}

\subsection{Compatibilities between products and coproducts}

We define another coproduct $\Delta_\haut$ on $\h$ in the following way: for all $x,y,z\in \h$,
$$\langle \Delta_\haut(x),y \otimes z\rangle=\langle x,z \haut y\rangle.$$

\begin{lemma}
For all forest $F \in \F$, $\displaystyle \Delta_\haut(F)=\sum_{\substack{F_1,F_2\in \F\\F_1F_2=F}} F_1 \otimes F_2$.
\end{lemma}

\begin{proof} Let $F,G,H \in \F$. Then:
\begin{eqnarray*}
\langle \Delta_\haut(F),f_G \otimes f_H\rangle&=&\langle F,f_H \haut f_G \rangle\\
&=&\langle F,f_{GH} \rangle\\
&=&\delta_{F,GH}\\
&=&\sum_{\substack{F_1,F_2\in \F\\F_1F_2=F}} \langle F_1 \otimes F_2, f_G \otimes f_H\rangle.
\end{eqnarray*}
As $(f_F)_{F\in \F}$ is a basis of $\h$ and $\langle-,-\rangle$ is non degenerate, this proves the result. \end{proof}\\

{\bf Remark.} As a consequence, the elements of $\T$ are primitive for this coproduct.\\

We now have defined three products, namely $m$, $\haut$, and $\bas$, and two coproducts, namely $\tdelta$ and $\tdelta_\haut$, on $\m$,
obtained from $\Delta$ and $\Delta_\haut$ by substracting their primitive parts. The following properties sum up the different compatibilities.

\begin{prop}
For all $x,y \in \m$:
\begin{eqnarray}
\label{E4} \tdelta(xy)&=&(x \otimes 1) \tdelta(y)+\tdelta(x)(1\otimes y)+x \otimes y,\\
\label{E5} \tdelta(x\haut y)&=&(x \otimes 1)\haut \tdelta(y)+\tdelta(x)\haut (1\otimes y)+x \otimes y,\\
\nonumber\\
\label{E6} \tdelta_\haut(xy)&=&(x \otimes 1) \tdelta_\haut(y)+\tdelta_\haut(x)(1\otimes y)+x \otimes y,\\
\label{E7} \tdelta_\haut(x \haut y)&=&(x \otimes 1) \haut \tdelta_\haut(y),\\
\label{E8} \tdelta_\haut(x \bas y)&=&(x \otimes 1) \bas \tdelta_\haut(y).
\end{eqnarray} \end{prop}

\begin{proof} It remains to consider the compatibility between $\haut$ or $\bas$  and $\tdelta_\haut$. Let $F,G\in \F-\{1\}$. We put $G=t_1\ldots t_n$,
where the $t_i$'s are trees. Then $F \haut G= (F \haut t_1) t_2\ldots t_n$, and $F \haut t_1$ is a tree. Hence:
\begin{eqnarray*}
\tdelta_\haut(F \haut G)&=&\sum_{i=1}^{n-1} (F \haut t_1)t_2 \ldots t_i \otimes t_{i+1}\ldots t_n\\
&=&\sum_{i=1}^{n-1} F \haut (t_1t_2 \ldots t_i) \otimes t_{i+1}\ldots t_n\\
&=&(F \otimes 1) \haut \tdelta_\haut(G).
\end{eqnarray*}
The proof is similar for $F \bas G$. So all these compatibilities are satisfied. \end{proof}\\

{\bf Remark.} There is no similar compatibility between $\tdelta$ and $\bas$. In particular, lemma \ref{19} is not available for $t\notin \T$.\\

This justifies the following definitions:

\begin{defi} \textnormal{\begin{enumerate}
\item A $\P_\haut$-bialgebra of type $1$ is a family $(A,m,\haut,\tdelta)$, such that:
\begin{enumerate}
\item $(A,m,\haut)$ is a $\P_\haut$-algebra.
\item $(A,\tdelta)$ is a coassociative, non counitary coalgebra.
\item Compatibilities (\ref{E4}) and (\ref{E5}) are satisfied.
\end{enumerate}
\item A $\P_\haut$-bialgebra of type $2$ is a family $(A,m,\haut,\tdelta_\haut)$, such that:
\begin{enumerate}
\item $(A,m,\haut)$ is a $\P_\haut$-algebra.
\item $(A,\tdelta_\haut)$ is a coassociative, non counitary coalgebra.
\item Compatibilities (\ref{E6}) and (\ref{E7}) are satisfied.
\end{enumerate}
\item A $\P_\bas$-bialgebra is a family $(A,m,\bas,\tdelta_\haut)$, such that:
\begin{enumerate}
\item $(A,m,\bas)$ is a $\P_\bas$-algebra.
\item $(A,\tdelta_\haut)$ is a coassociative, non counitary coalgebra.
\item Compatibilities (\ref{E6}) and (\ref{E8}) are satisfied.
\end{enumerate} \end{enumerate}} \end{defi}

{\bf Example.}  The augmentation ideal $\m$ of the infinitesimal Hopf algebra of trees $\h$ is both a $\P_\haut$-infinitesimal bialgebra 
of type $1$ and $2$, and also a $\P_\bas$-infinitesimal bialgebra.\\

If $A$ is a bialgebra of such a type, we denote by $Prim(A)$ the kernel of the coproduct. We deduce the definition of the following suboperads:

\begin{defi} \textnormal{
Let $n\in \mathbb{N}$. We put:
$$\left\{ \begin{array}{rcl}
\Q_\haut^{(1)}(n)&=&\left\{p\in \P_\haut(n)\:/\: \begin{array}{c}
\mbox{For all $A$, $\P_\haut $-infinitesimal bialgebra of type $1$,}\\
\mbox{and for $a_1,\ldots,a_n \in Prim(A)$,}\\
p.(a_1,\ldots,a_n)\in Prim(A).
\end{array}\right\},\\ \\
\Q_\haut^{(2)}(n)&=&\left\{p\in \P_\haut(n)\:/\: \begin{array}{c}
\mbox{For all $A$, $\P_\haut $-infinitesimal bialgebra of type $2$,}\\
\mbox{and for $a_1,\ldots,a_n \in Prim_\haut(A) $,}\\
p.(a_1,\ldots,a_n)\in Prim_\haut(A).
\end{array}\right\},\\ \\
\Q_\bas(n)&=&\left\{p\in \P_\bas(n)\:/\:\begin{array}{c}
\mbox{For all $A$, $\P_\bas $-infinitesimal bialgebra,}\\
\mbox{and for $a_1,\ldots,a_n \in Prim_\haut(A) $,}\\
p.(a_1,\ldots,a_n)\in Prim_\haut(A).
\end{array}\right\}. \end{array}\right.$$} \end{defi}

We identify $\P_\haut(n)$ and $\P_\bas(n)$ with the homogeneous component of weight $n$ of $\m$.
We put $Prim(\m)=Ker(\tdelta)$ and $Prim_\haut(\m)=Ker(\tdelta_\haut)$. We obtain:

\begin{prop} \label{27} \begin{enumerate}
\item For all $n \in \mathbb{N}$:
$$\Q_\haut^{(1)}(n)=\left\{p\in \P_\haut(n)\:/\: p\prodhaut(\tun,\ldots,\tun)\in Prim(\m)\right\}=\P_\haut(n)\cap Prim(\m).$$
\item For all $n \in \mathbb{N}$:
$$\Q_\haut^{(2)}(n)=\left\{p\in \P_\haut(n)\:/\: p\prodhaut(\tun,\ldots,\tun)\in Prim_\haut(\m)\right\}=\P_\haut(n)\cap Prim_\haut(\m).$$
\item For all $n \in \mathbb{N}$:
$$\Q_\bas(n)=\left\{p\in \P_\bas(n)\:/\: p\prodbas(\tun,\ldots,\tun)\in Prim_\haut(\m)\right\}=\P_\bas(n)\cap Prim_\haut(\m).$$
\end{enumerate} \end{prop}

\begin{proof} As $\m$ is a $\P_\haut$-infinitesimal bialgebra, by definition:
$$\Q_\haut^{(1)}(n)\subseteq \left\{p\in \P_\haut(n)\:/\: p\prodhaut(\tun,\ldots,\tun)\in Prim(\m)\right\}.$$
Moreover, $\left\{p\in \P_\haut(n)\:/\: p\prodhaut(\tun,\ldots,\tun)\in Prim(\m)\right\}=\P_\haut(n)\cap Prim(\m)$, as, for all $p\in \P_{\haut}(n)$,
$p\prodhaut(\tun,\ldots,\tun)=p \in \m$.

We now show that $\left\{p\in \P_\haut(n)\:/\: p\prodhaut(\tun,\ldots,\tun)\in Prim(\m)\right\}\subseteq \Q_\haut^{(1)}(n)$.
We take $p\in \P_\haut(n)$, such that $p\prodhaut(\tun,\ldots,\tun) \in Prim(\m)$. Let $\D=\{1,\ldots,n\}$ and let $A$ be the free $\P_\haut$-algebra
generated by $\D$ (with a unit). It can be described as the associative algebra $\h^\D$ generated by the set of planar rooted trees decorated by $\D$, 
and can be given a structure of $\P_\haut$-infinitesimal bialgebra. As $\m$ is freely generated by $\tun$ as a $\P_\haut$-algebra,
there exists a unique morphism of $\P_\haut$-algebras from $\m$ to $\m^\D$, augmentation ideal of $\h^\D$:
$$\xi:\left\{\begin{array}{rcl}
\m&\longrightarrow&\m^\D\\
\tun&\longrightarrow&\tdun{$1$}+\ldots+\tdun{$n$}.
\end{array}\right. $$
As $\tun \in Prim(\m)$ and $\tdun{$1$}+\ldots+\tdun{$n$} \in Prim(A)$, $\xi$ is a $\P_\haut$-infinitesimal bialgebra morphism from $\m$ to $\m^\D$.
So, $\xi(p\prodhaut(\tun,\ldots,\tun))\in Prim(A)$.

Let $F\in A$ be a forest, and $s_1 \geq_{h,l}\ldots \geq_{h,l} s_k$ its vertices. For all $i \in \{1,\ldots,k\}$, we put $d_i$ the decoration of $s_i$.
The {\it decoration word} associated to $F$ is the word $d_1\ldots d_n$. It belongs to $M(\D)$, the free monoid generated by the elements of $\D$.
For all $w \in M(\D)$, Let $A_w$ be the subspace of $A$ generated by forests whose decoration word is $w$. 
This defines a $M(\D)$-gradation of $A$, as a $\P_\haut$-infinitesimal bialgebra of type 1.

Consider the projection $\pi_{1,\ldots,n}$ onto $A_{1,\ldots,n}$. We get:
\begin{eqnarray*}
\pi_{1,\ldots,n} \circ \xi(p\prodhaut(\tun,\ldots,\tun))&\in &Prim(A),\\
&=&\pi_{1,\ldots,n}( p\prodhaut(\xi(\tun),\ldots,\xi(\tun)))\\
&=&\pi_{1,\ldots,n}( p\prodhaut(\tdun{$1$}+\ldots+\tdun{$n$},\ldots,
\tdun{$1$}+\ldots+\tdun{$n$}))\\
&=&p\prodhaut(\tdun{$1$},\ldots,\tdun{$n$}).
\end{eqnarray*}
So $p\prodhaut(\tdun{$1$},\ldots,\tdun{$n$})\in Prim(A)$.

Let $B$ be a $\P_\haut$-infinitesimal bialgebra and let $a_1,\ldots,a_n \in Prim(B)$. 
As $\m^\D$ is freely generated by the $\tdun{$i$}$'s, there exists a unique morphism of $\P_\haut$-algebras:
$$\chi:\left\{\begin{array}{rcl}
A&\longrightarrow&B\\
\tdun{$i$}&\longrightarrow&a_i.
\end{array}\right. $$
As the $\tdun{$i$}$ and the $a_i$'s are primitive, $\chi$ is a $\P$-infinitesimal bialgebra morphism. So:
$$\xi(p\prodhaut(\tdun{$1$},\ldots,\tdun{$n$}))=p.(\xi(\tdun{$1$}),\ldots,\xi(\tdun{$n$}))=p.(a_1,\ldots,a_n) \in \chi(prim(\m^\D))\subseteq Prim(A).$$
Hence, $p\in \Q_\haut^{(1)}(n)$.  The proof is similar for $\Q_\haut^{(2)}$ and $\Q_\bas$.  \end{proof}

\subsection{Suboperad $\Q_\haut^{(1)}$}

\begin{lemma} \label{28}
We define inductively the following elements of $\P_\haut$:
$$  \left\{\begin{array}{rcl}
q_1&=&\tun,\\
q_{n+1}&=&(\tun\tun-\tdeux)\circhaut(q_n,\tun)\:=\:q_n \tun-B^+(q_n), \mbox{ for } n \geq 1.
\end{array}\right. $$
Then, for all $n \geq 1$, $q_n$  belongs to $\Q_\haut^{(1)}$. Moreover, for all $x_1,\ldots,x_n \in Prim(\m)$:
$$\gamma(q_n\prodhaut(x_1,\ldots,x_n))=\gamma(x_1)x_2\ldots x_n.$$
\end{lemma}

{\bf Remark.}  These $q_n$'s are the same  as the $q_n$'s defined in section \ref{s3-2}. \\

\begin{proof} Let us remark that $f_{\tdeux}=\tun\tun-\tdeux \in Prim(\m)$. By proposition \ref{27}, $\tun\tun-\tdeux \in \Q_\haut^{(1)}(2)$.
As $\Q_\haut^{(1)}$ is a suboperad of $\P_\haut$, it follows that all the $q_n$'s belongs to $\Q_\haut^{(1)}(n)$.

Let $x_1,\ldots,x_n \in Prim(\m)$. Let us show that $\gamma(q_n\prodhaut(x_1,\ldots,x_n))=\gamma(x_1)x_2\ldots x_n$ by induction on $n$. 
If $n=1$, this is immediate. For $n=2$, $q_2\prodhaut(x_1x_2)=x_1x_2-x_1\haut x_2$. Moreover, $x_1\haut x_2$ is a linear span of forests 
whose first tree is not $\tun$. So $\gamma(q_2\prodhaut(x_1,x_2))=\gamma(x_1x_2)-0=\gamma(x_1)x_2$.\\

Suppose now the result true at rank $n-1$. Then:
\begin{eqnarray*}
q_n\prodhaut(x_1,\ldots,x_n)&=&q_2\prodhaut(\underbrace{q_{n-1}\prodhaut(x_1,\ldots,x_{n-1})}_{\in Prim(\m)},x_n),\\
\gamma(q_n\prodhaut(x_1,\ldots,x_n))&=&\gamma(q_2\prodhaut(q_{n-1}\prodhaut(x_1,\ldots,x_{n-1}),x_n))\\
&=&\gamma(q_{n-1}\prodhaut(x_1,\ldots,x_{n-1}))x_n\\
&=&\gamma(x_1)x_2\ldots x_n.
\end{eqnarray*}
\end{proof}

\begin{theo} \label{29}
The non-$\Sigma$-operad $\Q_\haut^{(1)}$ is freely generated by $\tdeux-\tun\tun$.
\end{theo}

\begin{proof} Let us first show that the family $(q_n)_{n\geq 1}$ generates $\Q_\haut^{(1)}$.  Let $\P$ be the suboperad of $\Q_\haut^{(1)}$
generated by the $q_n$'s. Let us prove by induction on $k$ that $\Q_\haut^{(1)}(k)=\P(k)$. If $k=1$, $\P(1)=\Q_\haut^{(1)}(1)=K \tun$. 
Suppose the result at all ranks $\leq k-1$. By the rigidity theorem for infinitesimal bialgebra of \cite{Loday2}, a basis of $\h$ is 
$(f_{t_1}\ldots f_{t_n})_{t_1\ldots t_n \in \F}$, so a basis of $Prim(\m)$ is:
$$\left( \gamma^{-1}_{Prim(\h)}(f_{t_1}\ldots f_{t_n})\right)_{t_1\ldots t_n \in \F}.$$
So, a basis of $\Q_\haut^{(1)}(k)$ is 
$\left( \gamma^{-1}_{Prim(\h)}(f_{t_1}\ldots f_{t_n})\right)_{\substack{t_1\ldots t_n \in \F\\ weight(t_1\ldots t_n)=k-1}}$.
By lemma \ref{28}:
$$\gamma^{-1}_{Prim(\h)}(f_{t_1}\ldots f_{t_n})=q_{n+1}\prodhaut(\tun,f_{t_1},\ldots f_{t_n}).$$
By the induction hypothesis, the $f_{t_i}$'s belongs to $\P$. So:
$$\gamma^{-1}_{Prim(\h)}(f_{t_1}\ldots f_{t_n})=q_{n+1}\circhaut(\tun, f_{t_1},\ldots f_{t_n})\in \P(n).$$
So $\Q_\haut^{(1)}=\P$. \\

Moreover, if we denote by $\P'$ the suboperad of $\Q_\haut^{(1)}$ generated by $q_2$, then, immediately, $\P'\subseteq \P$. 
Finally, by induction on $n$, $q_n \in \P'(n)$ for all $n\geq 1$ and $\P\subseteq \P'$. So $\P'=\P=\Q_\haut^{(1)}$ is generated by $q_2$.\\

Let $\P_{q_2}$ be the non-$\Sigma$-operad freely generated by $q_2$. There is a non-$\Sigma$-operad epimorphism:
$$ \Psi:\left\{\begin{array}{rcl}
\P_{q_2}&\longrightarrow &\Q_\haut^{(1)}\\
q_2&\longrightarrow&q_2.
\end{array}\right. $$
The dimension of $\P_{q_2}(n)$ is the number of planar binary rooted trees with $n$ leaves, that is to say the Catalan number 
$\displaystyle c_n=\frac{(2n-2)!}{(n-1)!n!}$. On the other side, the dimension of $\Q_\haut^{(1)}(n)$ is the number of planar rooted trees with $n$ vertices, 
that is to say $c_n$. So $\Psi$ is an isomorphism. \end{proof}\\

In other terms, in the language of \cite{LodayTriple}:

\begin{theo}
The triple of operads $(\mathbb{A}ss, \P_{\haut}^\Sigma, \mathbb{FREE}_2)$, where $\P_{\haut}^\Sigma$ is the symmetrisation of $\P_{\haut}$ 
and $\mathbb{FREE}_2$ is the free operad generated by an element in $\mathbb{FREE}_2(2)$, is a good triple of operads.
\end{theo}

{\bf Remark.} Note that if $A$ is a $\P_\haut$-bialgebra of type $1$, then $(A,m,\tdelta)$ is a non unitary infinitesimal bialgebra. 
Hence, if $(K\oplus A,m, \Delta)$ has an antipode $S$, then $-S$ is an eulerian idempotent for $A$.

\subsection{Another basis of $Prim(\h)$}

Recall that $\TT_b$ is freely generated (as a non-$\Sigma$-operad) by $\bdeux$. In particular, if $t_1,t_2 \in \T_b$, we denote:
$$t_1\vee t_2=\bdeux \circ (t_1,t_2).$$
Every element $t\in \T_b-\{\bun\}$ can be uniquely written as $t=t^l\vee t^r$.\\

There exists a morphism of operads:
$$\Theta : \left\{\begin{array}{rcl}
\T_b & \longrightarrow &\P_\haut\\
\bdeux & \longrightarrow &\tun\tun-\tdeux.
\end{array}\right. $$
By theorem \ref{29}, $\Theta$ is injective and its image is $\Q_\haut^{(1)}$. So, we obtain a basis of $\Q_\haut^{(1)}$ indexed by $\T_b$, 
given by $p_t=\Theta(t)$.  It is also a basis of $Prim(\m)$, which can be inductively computed by:
$$\left\{\begin{array}{rcl}
p_{\bun}&=&\tun,\\
p_{t_1\vee t_2}&=&(\tun \tun-\tdeux)\circhaut (p_{t_1},p_{t_2})
\:= \:p_{t_1}p_{t_2}-p_{t_1}\haut p_{t_2}.
\end{array}\right. $$

\subsection{From the basis $(f_t)_{t\in \T}$ to the basis $(p_t)_{t\in \T_b}$}

We define inductively the application $\kappa: \T_b \longrightarrow \T$ in the following way:
$$\kappa:\left\{\begin{array}{rcl}
\T_b&\longrightarrow&\T\\
\bun&\longrightarrow&\tun,\\
t_1\vee t_2&\longrightarrow & \kappa(t_2)\bas \kappa(t_1).
\end{array}\right. $$

{\bf Examples.}
$$\begin{array}{rclc|crclc|crclc|crcl}
\bdeux&\longrightarrow & \tdeux&&&\btroisun&\longrightarrow & \ttroisun&&&
\btroisdeux&\longrightarrow & \ttroisdeux&&&\bquatrecinq &\longrightarrow &\tquatredeux\\
\bquatreun &\longrightarrow &\tquatreun&&&\bquatredeux &\longrightarrow &\tquatretrois&&&
\bquatretrois &\longrightarrow &\tquatrequatre&&&\bquatrequatre &\longrightarrow &\tquatrecinq
\end{array}$$\\

It is easy to show that $\kappa$ is bijective, with inverse given by:
$$ \kappa^{-1}:\left\{\begin{array}{rcl}
\T&\longrightarrow&\T_b\\
\tun&\longrightarrow&\bun,\\
B^+(s_1\ldots s_m)&\longrightarrow&\kappa^{-1}(B^+(s_2\ldots s_m))\vee \kappa^{-1}(s_1).
\end{array}\right. $$

Let us recall the partial order $\leq$, defined in \cite{FoissyInf}, on the set $\F$ of planar forests, making it isomorphic to the Tamari poset.

\begin{defi} 
\textnormal{Let $F \in \F$. 
\begin{enumerate}
\item An {\it admissible transformation} on $F$ is a local transformation of $F$ of one of the following types
(the part of $F$ which is not in the frame remains unchanged):
$$\begin{array}{rccc}
\mbox{First kind: }&
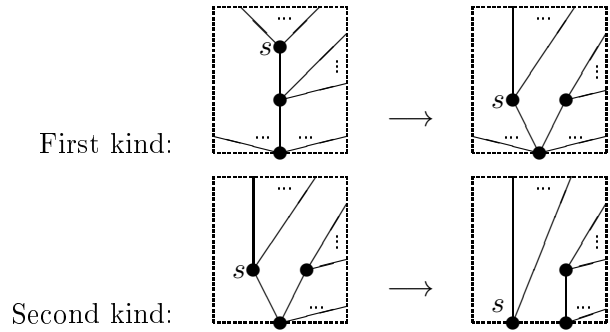
\begin{picture}(55,60)(-25,0)
\put(0,0){\circle*{5}} \put(0,0){\line(0,1){20}} \put(0,20){\circle*{5}} \put(0,20){\line(0,1){20}} \put(0,40){\circle*{5}} \put(0,40){\line(-1,1){15}}
\put(0,40){\line(1,1){15}} \put(-8,37){$s$} \put(-2,50){.} \put(0,50){.} \put(2,50){.} \put(0,0){\line(-4,1){25}} \put(0,0){\line(4,1){25}} \put(-10,5){.}
\put(-8,5){.} \put(-6,5){.} \put(10,5){.} \put(8,5){.} \put(6,5){.} \put(0,20){\line(4,1){25}} \put(0,20){\line(1,1){25}} \put(20,30){.} \put(20,32){.} \put(20,34){.}
\put(-25,0){\dashbox{1}(50,55)}\end{picture}
&\begin{picture}(22,0)(0,0)
\put(0,10){$ \longrightarrow $} 
\end{picture}
&\begin{picture}(55,60)(-25,0) 
\put(0,0){\circle*{5}} \put(0,0){\line(1,2){10}} \put(0,0){\line(-1,2){10}} \put(-10,20){\circle*{5}} \put(10,20){\circle*{5}} \put(-18,17){$s$} 
\put(0,0){\line(-4,1){25}} \put(0,0){\line(4,1){25}} \put(-14,5){.} \put(-12,5){.} \put(-10,5){.} \put(14,5){.} \put(12,5){.} \put(10,5){.} 
\put(-10,20){\line(0,1){35}} \put(-10,20){\line(2,3){23.5}} \put(-2,50){.} \put(0,50){.} \put(2,50){.} \put(10,20){\line(3,5){15}} \put(10,20){\line(4,1){15}}
\put(20,28){.} \put(20,30){.} \put(20,32){.} \put(-25,0){\dashbox{1}(50,55)} 
\end{picture}\\
\mbox{Second kind: }&
\begin{picture}(55,60)(-25,0)
\put(0,0){\circle*{5}} \put(0,0){\line(1,2){10}} \put(0,0){\line(-1,2){10}} \put(-10,20){\circle*{5}} \put(10,20){\circle*{5}} \put(-18,17){$s$} 
\put(0,0){\line(4,1){25}}\put(14,5){.} \put(12,5){.} \put(10,5){.} \put(-10,20){\line(0,1){35}} \put(-10,20){\line(2,3){23.5}} \put(-2,50){.} \put(0,50){.}
\put(2,50){.} \put(10,20){\line(3,5){15}} \put(10,20){\line(4,1){15}} \put(20,28){.} \put(20,30){.} \put(20,32){.} \put(-25,0){\dashbox{1}(50,55)} 
\end{picture}
&\begin{picture}(22,0)(0,0)
\put(0,10){$ \longrightarrow $}
\end{picture}
&\begin{picture}(55,60)(-25,0) \put(10,0){\circle*{5}} \put(10,20){\circle*{5}} \put(-10,0){\circle*{5}} \put(-18,4){$s$} \put(10,0){\line(0,1){20}}
\put(10,20){\line(3,5){15}} \put(10,20){\line(4,1){15}} \put(20,28){.} \put(20,30){.} \put(20,32){.} \put(-10,0){\line(0,1){55}} \put(-10,0){\line(2,5){22}}
\put(-2,50){.} \put(0,50){.} \put(2,50){.} \put(10,0){\line(3,1){15}} \put(-25,0){\dashbox{1}(50,55)} \put(13,5){.} \put(15,5){.} \put(17,5){.}
\end{picture}\end{array}$$
\item Let $F$ and $G \in \F$. We shall say that $F\leq G$ if there exists a finite sequence $F_0,\ldots,F_k$ of elements of $\F$ such that:
\begin{enumerate}
\item For all $i \in \{0,\ldots, k-1\}$, $F_{i+1}$ is obtained from $F_i$ by an admissible transformation.
\item $F_0=F$.
\item $F_k=G$.
\end{enumerate}\end{enumerate}}
\end{defi}

The aim of this section is to prove the following result:

\begin{theo} \label{32}
Let $t\in \T_b$. Then $\displaystyle p_t=\sum_{\substack{s\in \T \\ s\leq \kappa(t)}} f_s$.
\end{theo}

\begin{proof} By induction on the number $n$ of leaves of $t$. If $n=1$, then $t=\bun$ and $p_{\bun}=\tun=f_{\tun}$.
Suppose the result at all ranks $\leq n-1$. As $p_t$ is primitive, we can put:
$$ p_t=\sum_{s \in \T} a_s f_s.$$
Write $t=t_1 \vee t_2$. By the induction hypothesis:
$$p_{t_1}=\sum_{\substack{s_1\in \T \\ s_1\leq \kappa(t_1)}} f_{s_1}\mbox{ and }p_{t_2}=\sum_{\substack{s_2\in \T \\ s_2\leq \kappa(t_2)}} f_{s_2}.$$
As $t=t_1 \vee t_2$, $p_t=(\tun \tun -\tdeux)\circhaut (p_{t_1},p_{t_2})=p_{t_1}p_{t_2}-p_{t_1}\haut p_{t_2}$.
So, for all $s \in \T$, as $s$ is primitive for $\Delta_\haut$:
\begin{eqnarray*}
a_s&=&\langle p_t,s\rangle\\
&=&\langle p_{t_1}p_{t_2}-p_{t_1}\haut p_{t_2},s \rangle\\
&=&\langle p_{t_2} \otimes p_{t_1},\Delta(s)-\Delta_\haut(s)\rangle\\
&=&\langle p_{t_2} \otimes p_{t_1},\Delta(s)\rangle\\
&=&\sum_{\substack{s_1\in \T \\ s_1\leq \kappa(t_1)}} \sum_{\substack{s_2\in \T \\ s_2\leq \kappa(t_2)}}
\langle f_{s_2} \otimes f_{s_1} , \Delta(s)\rangle.
\end{eqnarray*}
So $a_s$ is the number of left-admissible cuts $c$ of $s$, such that $P^c(s) \leq \kappa(t_2)$ and $R^c(s) \leq \kappa(t_1)$.\\

Suppose that $a_s \neq 0$. Then, there exists a left-admissible cut $c$ of $s$, such that $P^c(s) \leq \kappa(t_2)$ and $R^c(s) \leq \kappa(t_1)$. 
As $s$ is a tree, $s \leq \kappa(t_2)\bas \kappa(t_1)=\kappa(t)$. Moreover, by considering the degree of $P^c(s)$, this cut $c$ is unique, so $a_s=1$.
Reciproquely, if $s \leq \kappa(t)$, if $c$ is the unique left admissible cut such that $weight(P^c(s))=weight(t_2)$,  then $P^c(s)\leq \kappa(t_2)$ 
and $R^c(s) \leq \kappa(t_1)$. So $a_s \neq 0$. Hence, $(s \leq \kappa(t)) \Longrightarrow (a_s \neq 0) \Longrightarrow (a_s=1) 
\Longrightarrow (s \leq \kappa(t))$.  This proves theorem \ref{32}. \end{proof}\\

Let $\mu$ be the Möbius function of the poset $\F$ (\cite{Stanley,Stanley2}). By the Möbius inversion formula:

\begin{cor}
Let $s \in \T$. Then $\displaystyle f_s=\sum_{t\in T_b,\: \kappa(t)\leq s} \mu(\kappa(t),s) p_t$.
\end{cor}

\subsection{Suboperad $\Q_\haut^{(2)}$}

For all $n \in \mathbb{N}$, we put $c_{n+1}=B^+(\tun^n)$. In other terms, $c_{n+1}$ is the corolla tree with $n+1$ vertices, or equivalently with $n$ leaves.\\

{\bf Examples.} $c_1=\tun$, $c_2=\tdeux$, $c_3=\ttroisun$, $c_4=\tquatreun$, $c_5=\tcinqun$\ldots

\begin{lemma} \label{34}
The set $\T$ is a basis of the operad $\Q_\haut^{(2)}$. As an operad, $\Q_\haut^{(2)}$ is generated by the $c_n$'s, $n\geq 2$.
Moreover, for all $k,l \geq 2$, 
$$c_k \circhaut (c_l,\underbrace{\tun,\ldots,\tun}_{\mbox{\scriptsize $k-1$ times}})
=c_l \circhaut (\underbrace{\tun,\ldots,\tun}_{\mbox{\scriptsize $l-1$ times}},c_k).$$
\end{lemma}

\begin{proof} The operad $\Q_\haut^{(2)}$ is identified with $Prim_\haut(\m)$ by proposition \ref{27}. So $Prim_\haut(\m)$ is equal to $Vect(\T)$.
Let $\P$ be the suboperad of $\Q_\haut^{(2)}$ generated by the corollas. Let $t \in \T$, of weight $n$. Let us prove that $t \in \P$ by induction on $n$.
If $n=1$, then $t=\tun \in \P$. If $n\geq 2$, we can suppose that $t=B^+(t_1\ldots t_k)$, with $t_1,\ldots,t_k \in \P$. Then, by theorem \ref{11}:
$$c_{k+1} \circhaut(t_1,\ldots,t_k,\tun)=(\tun^k\circhaut(t_1,\ldots,t_k))\haut \tun=(t_1\ldots t_k) \haut \tun=B^+(t_1\ldots t_k)=t.$$
So $t \in \P$. hence, $\P=\Q_\haut^{(2)}$.\\

Let $k,l\geq 2$. Then, by theorem \ref{11}:
\begin{eqnarray*}
c_k \circhaut(c_l,\tun,\ldots,\tun)&=&(\tun^{k-1}\circhaut(c_l,\tun,\ldots,\tun)) \haut \tun\\
&=&(c_l \tun^{k-2})\haut \tun\\
&=&B^+(c_l \tun^{k-2})\\
&=&B^+(B^+(\tun^{l-1}) \tun^{k-2}).
\end{eqnarray*}
On the other hand:
\begin{eqnarray*}
c_l \circhaut(\tun,\ldots,\tun,c_k)&=&
(\tun^{l-1}\circhaut(\tun,\ldots,\tun)) \haut c_k\\
&=&(\tun^{l-1})\haut c_k\\
&=&(\tun^{l-1})\haut B^+(c^{k-1})\\
&=&B^+(((\tun^{l-1})\haut \tun) \tun^{k-2})\\
&=&B^+(B^+(\tun^{l-1}) \tun^{k-2}).
\end{eqnarray*}
So $c_k \circhaut(c_l,\tun,\ldots,\tun)=c_l \circhaut(\tun,\ldots,\tun,c_k)$. \end{proof}

\begin{defi} \textnormal{
The operad $\TT$ is the non-$\Sigma$-operad generated by elements $c_n \in \TT(n)$, for $n \geq 2$, 
and the following relations: for all $k,l \geq 2$,
$$c_k \circ (c_l,\underbrace{I,\ldots,I}_{\mbox{\scriptsize $k-1$ times}})=c_l \circ (\underbrace{I,\ldots,I}_{\mbox{\scriptsize $l-1$ times}},c_k).$$
In other terms, a $\TT$-algebra $A$ has a family of $n$-multilinear products $[.,\ldots,.]:A^{\otimes n} \longrightarrow A$ for all $n \geq 2$, 
with the associativity condition:
$$[[a_1,\ldots,a_l],a_{l+1},\ldots,a_{l+k}]=[a_1,\ldots,a_{l-1},[a_l,\ldots,a_{l+k}]].$$
In particular, $[.,.]$ is associative.}
\end{defi}

\begin{theo}
The operads $\TT$ and $\Q_\haut^{(2)}$ are isomorphic.
\end{theo}

\begin{proof} By lemma \ref{34}, there is an epimorphism of operads:
$$  \left\{ \begin{array}{rcl}
\TT&\longrightarrow & \Q_\haut^{(2)} \\
c_n&\longrightarrow &c_n.
\end{array} \right.$$
In order to prove this is an isomorphism, it is enough to prove that $dim(\TT(n))\leq dim(\Q_\haut^{(2)}(n))$ for all $n \geq 2$. 
By lemma \ref{34}, $dim(\Q_\haut^{(2)}(n))$ is the $n$-th Catalan number. Because of the defining relations, $\TT(n)$ is generated 
as a vector space by elements of the form $c_l\circ(I,b_2,\ldots,b_l)$, with $b_i \in \TT(n_i)$,  such that $n_1+\ldots+n_l=n-1$. 
Hence, we define inductively the following subsets the free non-$\Sigma$-operad generated by the $c_n$'s, $n\geq 2$:
$$X(n)= \left\{ \begin{array}{l}
\{I\}\mbox{ if } n=1,\\
\displaystyle \: \bigcup_{l=2}^{n} \: \bigcup_{i_2+\ldots+i_l=n-1}c_l \circ (I,X(i_2),\ldots,X(i_l) )\mbox{ if } n\geq 2.
\end{array}\right.$$
Then the images of the elements of $X(n)$ linearly generate $\TT(n)$, so $dim(\TT(n))\leq card(X(n))$ for all $n$. 
We now put $a_n=card(X(n))$ and prove that $a_n$ is the $n$-th Catalan number. We denote by $A(h)$ their generating formal series. Then:
$$  \left\{ \begin{array}{l}
a_1=1,\\
\displaystyle a_n=\sum_{l=2}^{n}\: \sum_{i_2+\ldots+i_l=n-1} a_{i_1}\ldots a_{i_l} \mbox{ if } n\geq 2.
\end{array}\right.$$
In terms of generating series:
$$A(h)-a_1h=h \frac{A(x)}{1-A(x)}.$$
So $A(h)^2-A(h)+h=0$. As $A(h)=1$:
$$A(h)=\frac{1-\sqrt{1-4h}}{2}.$$
So $a_n$ is the $n$-th Catalan number for all $n \geq 1$. \end{proof}\\

In other terms:

\begin{theo}
The triple of operads $(\mathbb{A}ss, \P_{\haut}^\Sigma, \TT^\Sigma)$ is a good triple of operads.
\end{theo}

{\bf Remark.} Note that if $A$ is a $\P_\haut$-bialgebra of type $2$, then $(A,m,\tdelta_\haut)$ is a non unitary infinitesimal bialgebra. 
Hence, if $(K\oplus A,m, \Delta_\haut)$ has an antipode $S_\haut$, then $-S_\haut$ is an eulerian idempotent for $A$.

\subsection{Suboperad $\Q_\bas$}

\begin{lemma}
The set $\T$ is a basis of the operad $\Q_\bas$. As an operad, $\Q_\bas$ is generated by $\tdeux$.
\end{lemma}

\begin{proof} Let $\P$ be the suboperad of $\Q_\bas$ generated by $\tdeux$. Let $t \in \T$, of weight $n$. Let us prove that $t \in \P$
by induction on $n$.  If $n=1$ or $2$, this is obvious. If $n \geq 2$, suppose that $t=B^+(t_1\ldots t_k)$. By the induction hypothesis, 
$t_1$ and $B^+(t_2\ldots t_k)$ belong to $\P$. Then:
$$t=t_1\bas B^+(t_2\ldots t_k)=\tdeux \circbas (t_1,B^+(t_2\ldots t_k)).$$
So $t\in \P$. \end{proof}

\begin{theo}
The non-$\Sigma$-operad $\Q_\bas$ is freely generated by $\tdeux$.
\end{theo}

\begin{proof} Similar as the proof of theorem \ref{29}. \end{proof}\\

In other terms:

\begin{theo}
The triple of operads $(\mathbb{A}ss, \P_{\bas}^\Sigma, \mathbb{F}_2)$, where $\mathbb{F}_2$ is the free operad generated by an element 
in $\mathbb{F}_2(2)$, is a good triple of operads.
\end{theo}

{\bf Remark.} Note that if $A$ is a $\P_\bas$-bialgebra, then $(A,m,\tdelta)$ is a non unitary infinitesimal bialgebra. 
Hence, if $(K\oplus A,m, \Delta)$ has an antipode $S$, then $-S$ is an eulerian idempotent for $A$.

\section{A rigidity theorem for $\P_\haut$-algebras}

\subsection{Double $\P_\haut$-infinitesimal bialgebras}

\begin{defi} \textnormal{
A double $\P_\haut$-infinitesimal bialgebra is a family $(A,m,\haut,\tdelta,\tdelta_\haut)$ where
$m,\haut:A\otimes A\longrightarrow A$, $\tdelta,\tdelta_\haut:A\longrightarrow A\otimes A$, with the following compatibilities:
\begin{enumerate}
\item  $(A,m,\haut)$ is a (non unitary) $\P_\haut$-algebra.
\item For all $x \in A$:
$$  \left\{ \begin{array}{rcl}
(\tdelta \otimes Id) \circ \tdelta(x)&=&(Id \otimes \tdelta) \circ \tdelta(x),\\
(\tdelta_\haut \otimes Id) \circ \tdelta_\haut(x)&=&(Id \otimes \tdelta_\haut) \circ \tdelta_\haut(x),\\
(\tdelta \otimes Id) \circ \tdelta_\haut(x)&=&(Id \otimes \tdelta_\haut) \circ \tdelta(x).
\end{array}\right.$$
In other terms, $(A,\tdelta^{cop},\tdelta_\haut^{cop})$ is a $\P_\haut$-coalgebra.
\item $(A,m, \haut,\tdelta)$ is a $\P_\haut$-bialgebra of type $1$.
\item $(A,m, \haut,\tdelta_\haut)$ is a $\P_\haut$-bialgebra of type $2$.
\end{enumerate}} \end{defi}

{\bf Remark.} If $(A,m,\haut,\tdelta,\tdelta_\haut)$ is a graded double $\P_\haut$-infinitesimal bialgebra, 
with finite-dimensional homogeneous components, then its graded dual $(A^*,\tdelta^{*,op}, \tdelta_\haut^{*,op},m^{*,cop},\haut^{*,cop})$ also is.

\begin{theo}
$(\m,m,\haut,\tdelta,\tdelta_\haut)$ is a double $\P$-infinitesimal bialgebra.
\end{theo}

\begin{proof} We already now that $(\m,m,\haut)$ is a $\P_\haut$-algebra. Moreover, $(\m,\tdelta^{cop},\tdelta_\haut^{cop})$ is isomorphic to
 $(\m^*,m^*,\haut^*)$  via the pairing $\langle-,-\rangle$, so it is a $\P_\haut$-coalgebra. It is already known that $(\m,m,\tdelta)$ and $(\m,\haut,\tdelta)$ are infinitesimal bialgebras. As $(\m,\haut,\tdelta)$ is isomorphic to $(\m^{op},m^{op},\tdelta_\haut^{cop})$ via the pairing
$\langle-,-\rangle$, it is also an infinitesimal bialgebra. So all the needed compatibilities are satisfied. \end{proof}\\

{\bf Remarks.} \begin{enumerate}
\item  Via the pairing $\langle-,-\rangle$, $\m$ is isomorphic to its graded dual as an double $\P_\haut$-infinitesimal bialgebra.
As a consequence, as $\m$ is the free $\P_\haut$-algebra generated by $\tun$, 
then $\m^{cop}$ is also the cofree $\P_\haut$-coalgebra cogenerated by $\tun$.
\item All these results can be easily extended to infinitesimal Hopf algebras of decorated planar rooted trees, 
in other terms to every free $\P_\haut$-algebras.
\end{enumerate}

\begin{lemma} \label{43} 
In the double infinitesimal $\P_\haut$-algebra $\m$, $Ker(\tdelta)\cap Ker(\tdelta_\haut)=Vect(\tun)$.
\end{lemma}

\begin{proof} $\supseteq$. Obvious. 

$\subseteq$. Let $x \in Ker(\tdelta)\cap Ker(\tdelta_\haut)$. Then $\tdelta_\haut(x)=0$, so $x$ is a linear span of trees. We can write:
$$x=\sum_{t\in \T} a_t t.$$
Consider the terms in $\m \otimes \tun$ of $\tdelta(x)$. We get $\displaystyle \sum_{t\in \T-\{\tun\}} a_t B^-(t) \otimes \tun=0$,
where $B^-(t)$ is the forest obtained by deleting the root of $t$. So, if $t \neq \tun$, then $a_t=0$. So $x \in vect(\tun)$. \end{proof}\\

{\bf Remark.} This lemma can be extended to any free $\P_\haut$-algebra: if $V$ is a vector space, 
then the free $\P_\haut$-algebra $F_{\P_\haut}(V)$ generated by $V$ is given a structure of double $\P_\haut$-infinitesimal bialgebra 
by $\tdelta(v)=\tdelta_\haut(v)=0$ for all $v \in V$. In this case, $Ker(\tdelta) \cap Ker(\tdelta_\haut)=V$ for $F_{\P_\haut}(V)$.

\subsection{Connected double $\P_\haut$-infinitesimal bialgebras}

{\bf Notations.} Let $A$ be a double $\P_\haut$-infinitesimal bialgebra. The iterated coproducts will be denoted in the following way: 
for all $n \in \mathbb{N}$,
\begin{eqnarray*}
\tdelta^n &:& \left\{ \begin{array}{rcl}
A&\longrightarrow & A^{\otimes(n+1)}\\ 
a&\longrightarrow &a^{(1)}\otimes \ldots \otimes a^{(n+1)},
\end{array}\right.\\ \\
\tdelta_\haut^n &:& \left\{ \begin{array}{rcl}
A&\longrightarrow & A^{\otimes(n+1)}\\ 
a&\longrightarrow &a_\haut^{(1)}\otimes \ldots \otimes a_\haut^{(n+1)}.
\end{array}\right.
\end{eqnarray*}

\begin{defi}
\textnormal{Let $A$ be a double $\P_\haut$-infinitesimal bialgebra. It will be said {\it connected} if, for any $a\in A$, 
every iterated coproduct $A\longrightarrow A^{\otimes (n+1)}$ vanishes on $a$ for a great enough $n$.}
\end{defi}

\begin{theo}
Let $A$ be a connected double $\P_\haut$-infinitesimal bialgebra. Then $A$ is isomorphic to the free $\P_\haut$-algebra 
generated by $Prim(A)=Ker(\tdelta)\cap Ker(\tdelta_\haut)$ as a double $\P_\haut$-infinitesimal bialgebra.
\end{theo}

\begin{proof} {\it First step.} We shall use the results on infinitesimal Hopf algebras of \cite{FoissyInf}. We show that $A=Prim(A)+A.A+A\haut A$.
As $(A,\haut,\tdelta)$ is a connected non unitary infinitesimal bialgebra, it (or more precisely its unitarisation) has an antipode $S_\haut$, defined by:
$$S_\haut : \left\{ \begin{array}{rcl}
A&\longrightarrow & A\\
a&\longrightarrow & \displaystyle \sum_{i=0}^{\infty} (-1)^{i+1} a^{(1)} \haut \ldots \haut a^{(i+1)}.
\end{array}\right.$$  
As $(A,\tdelta)$ is connected, this makes sense. Moreover, $-S_\haut$ is the projector on $Ker(\tdelta)$ in the direct sum 
$A=Ker(\tdelta)\oplus A \haut A$.

In the same order of idea, as $(A,m,\tdelta_\haut)$ is a connected infinitesimal bialgebra, we can define its antipode $S^\haut$ by:
$$S^\haut: \left\{ \begin{array}{rcl}
A&\longrightarrow & A\\
a&\longrightarrow & \displaystyle \sum_{i=0}^{\infty} (-1)^{i+1} a^{(1)}_\haut \ldots a^{(i+1)}_\haut,
\end{array}\right.$$ 
and $-S^\haut$ is the projector on $Ker(\tdelta_\haut)$ in the direct sum $A=Ker(\tdelta_\haut)\oplus A.A$.

Let $a\in A$, $b \in Ker(\tdelta_\haut)$. Then $\tdelta_\haut(a \haut b)=(a \otimes 1)\tdelta_\haut(b)=0$.
So $A \haut Ker(\tdelta_\haut)$ is a subset of  $Ker(\tdelta_\haut)$. Moreover, if $\tdelta_\haut(a)=0$, 
then $(Id \otimes \tdelta_\haut) \circ \tdelta(a)=(\tdelta \otimes Id) \circ \tdelta_\haut(a)=0$.
So $\tdelta(a) \in A \otimes Ker(\tdelta_\haut)$. As a consequence, if $n \geq 1$:
$$\tdelta^n(a)=(\tdelta^{n-1} \otimes Id) \circ \tdelta(a) \in A^{\otimes n} \otimes Ker(\tdelta_\haut).$$
Hence, for all $n \in \mathbb{N}$, $\tdelta^n(Ker(\tdelta_\haut)) \in A^{\otimes n} \otimes Ker(\tdelta_\haut)$.
Finally, we deduce that $S_\haut(Ker(\tdelta_\haut)) \subseteq Ker(\tdelta_\haut)$.

Let $a \in A$. Then $S^\haut(a) \in Ker(\tdelta_\haut)$ and $S_\haut \circ S^\haut(a) \in Ker(\tdelta) \cap Ker(\tdelta_\haut)=Prim(A)$ 
by the preceding point. Moreover:
\begin{eqnarray*}
S^\haut(a)&=&-a+A.A,\\
S_\haut\circ S^\haut(a)&=&-S^\haut(a)+A\haut A,\\
S_\haut \circ S^\haut(a)&=&a+A.A+A\haut A.
\end{eqnarray*}
This proves the first step. \\

{\it Second step.} As $A$ is connected, it classically inherits a filtration of $\P_\haut$-algebra given by the kernels of the iterated coproducts. 
We denote by $deg_p$ the associated degree. In particular, for all $x \in A$, $deg_p(x)\leq 1$ if, and only if, $x\in Prim(A)$.
Let $B$ be the $\P_\haut$-subalgebra of $A$ generated by $Prim(A)$. Let $a\in A$, let us show that $a \in B$ by induction on $n=deg_p(a)$. 
If $n\leq 1$, then $a \in Prim(A)\subseteq B$. Suppose the result true at all ranks $\leq n-1$. Then, by the first step, we can write:
$$a=b+\sum_i a_i b_i+\sum_j c_jd_j,$$
with $b \in Prim(A)$, $a_i$, $b_i$, $c_j$, $d_j \in A$. Because of the filtration, we can suppose that $deg_p(a_i)$, $deg_p(b_i)$, 
$deg_p(c_j)$, $deg_p(d_j)<n$. By the induction hypothesis, they belong to $B$, so $a\in B$. \\

{\it Last step.} So, there is an epimorphism of $\P_\haut$-algebras:
$$ \phi : \left\{ \begin{array}{rcl}
F_{\P_\haut}(Prim(A))&\longrightarrow & A\\
a \in Prim(A) &\longrightarrow &a,
\end{array}\right.$$
where $F_{\P_\haut}(Prim(A))$ is the free $\P_\haut$-algebra generated by $Prim(A)$. 
As the elements of $Prim(A)$ are primitive both in $A$ and $F_{\P_\haut}(Prim(A))$, this is a morphism of double $\P_\haut$-infinitesimal
bialgebras. Suppose that it is not monic. Take then $x \in Ker(\phi)$, non-zero, of minimal degree.
Then it is primitive, so belongs to $Prim(A)$ (lemma \ref{43}). Hence, $\phi(a)=a=0$: this is a contradiction. So $\phi$ is a bijection. \end{proof}\\

In other terms:

\begin{cor}
The triple of operads $\left((\P_\haut^\Sigma)^{op},\P_\haut^\Sigma,\Vect\right)$ is a good triple. Here, $\Vect$ denotes the operad of vector spaces:
$$\Vect(k)=\left\{ \begin{array}{rcl}
KI&\mbox{if} &k=1,\\
0&\mbox{if}&k\neq 1,
\end{array}\right.$$
where $I$ is the unit of $\Vect$.
\end{cor}

We also showed that $S_\haut \circ S^\haut$ is the projection on $Prim(A)$ in the direct sum $A=Prim(A)\oplus(A.A+A\haut A)$. 
So $S_\haut \circ S^\haut$ is the eulerian idempotent.

\bibliographystyle{amsplain}
\bibliography{biblio}

\providecommand{\bysame}{\leavevmode\hbox to3em{\hrulefill}\thinspace}
\providecommand{\MR}{\relax\ifhmode\unskip\space\fi MR }
\providecommand{\MRhref}[2]{%
  \href{http://www.ams.org/mathscinet-getitem?mr=#1}{#2}
}
\providecommand{\href}[2]{#2}
\begin{thebibliography}{10}

\bibitem{Connes}
Alain Connes and Dirk Kreimer, \emph{Hopf algebras, {R}enormalization and
  {N}oncommutative geometry}, Comm. Math. Phys \textbf{199} (1998), no.~1,
  203--242, arXiv:hep-th/98 08042.

\bibitem{FoissyOperades}
Lo{\"\i}c Foissy, \emph{Koszularity of the operads of forests}, in preparation.

\bibitem{Foissy2}
\bysame, \emph{Les alg\`ebres de {H}opf des arbres enracin\'es, {I}}, Bull.
  Sci. Math. \textbf{126} (2002), 193--239.

\bibitem{FoissyInf}
\bysame, \emph{The infinitesimal {H}opf algebra and the poset of planar rooted
  forests}, J. Algebraic Combin. (2009), arXiv: 08 02.0442.

\bibitem{Holtkamp}
Ralf Holtkamp, \emph{Comparison of {H}opf {A}lgebras on {T}rees}, Arch. Math.
  (Basel) \textbf{80} (2003), no.~4, 368--383.

\bibitem{Kreimer1}
Dirk Kreimer, \emph{On the {H}opf algebra structure of pertubative quantum
  field theories}, Adv. Theor. Math. Phys. \textbf{2} (1998), no.~2, 303--334,
  arXiv:q-alg/97 07029.

\bibitem{Kreimer2}
\bysame, \emph{On {O}verlapping {D}ivergences}, Comm. Math. Phys. \textbf{204}
  (1999), no.~3, 669--689, arXiv:hep-th/98 10022.

\bibitem{Kreimer3}
\bysame, \emph{Combinatorics of (pertubative) {Q}uantum {F}ield {T}heory},
  Phys. Rep. \textbf{4--6} (2002), 387--424, arXiv:hep-th/00 10059.

\bibitem{LodayTriple}
Jean-Louis Loday, \emph{Generalized bialgebras and triples of operads}, arXiv:
  math/061 1885, 2006.

\bibitem{Loday2}
Jean-Louis Loday and Maria~O. Ronco, \emph{On the structure of cofree hopf
  algebras}, J. Reine Angew. Math. \textbf{592} (2006), 123--155, available at
  http://www-irma.u-strasbg.fr/$\tilde{\:}$loday/.

\bibitem{Markl}
Martin Markl, Steve Shnider, and Jim Stasheff, \emph{Operads in algebra,
  topology and physics}, Mathematical Surveys and Monographs, vol.~96, American
  Mathematical Society, 2002.

\bibitem{Stanley}
Richard~P. Stanley, \emph{Enumerative combinatorics. {V}ol. 1.}, Cambridge
  Studies in Advanced Mathematics, no.~49, Cambridge University Press,
  Cambridge, 1997.

\bibitem{Stanley2}
\bysame, \emph{Enumerative combinatorics. {V}ol. 2.}, Cambridge Studies in
  Advanced Mathematics, no.~62, Cambridge University Press, Cambridge, 1999.

\end{thebibliography}

\end{document}